\documentclass[11pt]{article}

\usepackage{geometry}
\usepackage[utf8]{inputenc}
\usepackage{amsthm, amssymb, natbib, graphicx}
\usepackage{listings}
\usepackage[ruled,vlined]{algorithm2e}
\usepackage{amsmath,tikz}
\usepackage{array, makecell}
\usepackage{sidecap, caption}
\usepackage{xcolor, color} 
\usepackage{todonotes}
\usepackage{comment} 
\usepackage{algorithm2e}
\usepackage[most]{tcolorbox}
\usepackage{lipsum}
\usepackage[compact]{titlesec}
\setcitestyle{numbers,comma,square}
\usepackage{caption}
\usepackage{tabularray}
\usepackage{diagbox}
\UseTblrLibrary{booktabs}
\usepackage{soul}
\usepackage{longtable}

\newtheorem{theorem}{Theorem}[section]
\newtheorem{corollary}{Corollary}[section]
\newtheorem{lemma}[theorem]{Lemma}
\newtheorem{example}[theorem]{Example}
\newtheorem{remark}{Remark}[section]

\newtheorem{definition}{Definition}[section]
\DeclareMathOperator*{\argmin}{argmin}

\newcommand{\R}{\mathbb{R}}

\newcommand{\rr}[1]{\textcolor{red}{#1}}

\usepackage[normalem]{ulem}
\newcommand{\bb}[1]{\textcolor{blue}{#1}}

\newcommand{\ignore}[1]{}

\title{Second-order methods for quartically-regularised cubic polynomials, with applications to high-order tensor methods
}

\author{\thanks{wenqi.zhu@maths.ox.ac.uk, University of Oxford, UK} Wenqi Zhu,  \thanks{cartis@maths.ox.ac.uk, University of Oxford, UK} Coralia Cartis}

\author{Coralia Cartis\thanks{The order of the authors is alphabetical; the second author (Wenqi Zhu) is the primary contributor.} \textsuperscript{\normalfont ,}\thanks{Mathematical Institute, Woodstock Road, University of Oxford, Oxford, UK, OX2 6GG.  {\tt cartis@maths.ox.ac.uk} } \quad and \quad
Wenqi Zhu\footnotemark[1] \textsuperscript{\normalfont,}\thanks{Mathematical Institute, Woodstock Road, University of Oxford, Oxford, UK, OX2 6GG.
{\tt wenqi.zhu@maths.ox.ac.uk}}}
\date{August 30, 2023; revised December 10, 2024}

\begin{document}

\maketitle

\begin{abstract}
There has been growing interest in high-order tensor methods for nonconvex optimization, with adaptive regularization, as they possess better/optimal worst-case evaluation complexity globally and faster convergence asymptotically. These algorithms crucially rely on repeatedly minimizing nonconvex multivariate Taylor-based polynomial sub-problems, at least locally. Finding efficient techniques for the solution of these sub-problems, beyond the second-order case, has been an open question.
This paper proposes a second-order method, Quadratic Quartic Regularisation (QQR), for efficiently minimizing nonconvex quartically-regularized cubic polynomials, such as the AR$p$ sub-problem \cite{birgin2017worst} with $p=3$. 
Inspired by \cite{Nesterov2022quartic},  QQR approximates the third-order tensor term by a linear combination of quadratic and quartic terms, yielding (possibly nonconvex) local models that are solvable to global optimality. In order to achieve accuracy $\epsilon$ in the {first-order criticality of the sub-problem in finitely many iterations}, we show that the error in the QQR method decreases either linearly or by at least $\mathcal{O}(\epsilon^{4/3})$ for locally convex iterations, while in the {nonconvex case}, by at least $\mathcal{O}(\epsilon)$; thus improving, on these types of iterations, the general cubic-regularization bound. Preliminary numerical experiments indicate that two QQR variants perform competitively with state-of-the-art approaches such as ARC  (also known as AR$p$ with $p=2$), achieving either lower objective value or iteration counts.
\end{abstract}

\section{Introduction and Problem Set-up}
In this paper, we consider the unconstrained nonconvex optimization problem, 
\begin{equation}
  \min_{x\in \R^n} f(x)  
  \label{min f}
\end{equation}
where  $f: \R^n \rightarrow \R$ is $p$-times continuously differentiable and bounded below, $p \ge 1$. 
Recent works~\cite{birgin2017worst, cartis2020concise, cartis2020sharp, Nesterov2021implementable} have shown that some optimization algorithms have better worst-case evaluation complexity when using high-order derivative information of the objective function, coupled with regularization techniques. In these methods, typically, a  polynomial local model $m_p(x, s)$ is constructed that approximates $f(x+s)$ {around} the iterate $x=x_k$. Then, $x_k$ is iteratively updated by $
    s_k \approx \argmin_{s\in \R^n}{m_p(x_k, s)}, \text{ and } x_{k+1} := x_k+s_k, 
$
whenever sufficient objective decrease is obtained, until an approximate local minimiser of $f$ is generated, satisfying $\|\nabla_x f(x_k)\|\leq \epsilon_g$ and $\lambda_{\min}\left(\nabla^2_x f(x_k)\right)\geq -\epsilon_H$. The model $m_p$ is based on the $p$th-order Taylor expansion of $f(x_k+s)$ at $x_k${, that is} $
   T_p(x_k, s) := f(x_k) + \sum_{j=1}^p \frac{1}{j!} \nabla_x^j f(x_k)[s]^j,
$
where $\nabla^j f(x_k) \in \R^{n^j}$ is a $j$th-order tensor and $\nabla^j f(x_k)[s]^j$ is the $j$th-order derivative of $f$ at $x_k$ along $s \in \R^n$. To ensure a bounded below local model\footnote{Unless otherwise stated, $\|\cdot\|$ denotes the Euclidean norm in this paper.}, and a convergent ensuing method,  an adaptive regularization term is added to $T_p$ leading to 
\begin{equation}
    m_p(x_k, s) = T_p(x_k, s)+ \frac{\sigma_k}{p+1}\|s\|^{p+1},
    \label{subprob}
\end{equation}
where $\sigma_k > 0$. The case of $p=1$ in \eqref{subprob} gives a steepest descent model and the case of $p=2$ gives a Newton-like model; in this paper, we focus on the case of $p =3$. The construction \eqref{subprob} corresponds to the sub-problem in the well-known adaptive regularization algorithmic framework AR$p$~\cite{birgin2017worst, cartis2020sharp, cartis2020concise}. In AR$p$, the parameter $\sigma_k$ is adjusted adaptively to ensure progress towards optimality over the iterations. 
Under Lipschitz continuity assumptions on $\nabla^p f(x)$, AR$p$ requires no more than $\mathcal{O}\left({\max\left[\epsilon_g^{-\frac{p+1}{p}},\epsilon_H^{-\frac{p+1}{p-1}} \right]}\right)$ evaluations of $f$ and its derivatives to compute a local minimizer to accuracy $(\epsilon_g, \epsilon_H)$ for  first order and second order criticality. This result shows that to achieve the same criticality conditions,  the evaluation complexity bound improves as we increase the order $p$. For instance, the AR$3$ algorithm exhibits better worst-case performance than first/second order methods, with evaluation complexity\footnote{The complexity excludes the computational cost of minimizing the sub-problem \eqref{subprob}.} of $\mathcal{O}\left(\max\{\epsilon_g^{-4/3}, \epsilon_H^{-2}\}\right)$. 
{This superior complexity result (compared to first- and second-order methods) motivates us to find efficient ways to solve subproblem \eqref{subprob} for $p=3$, which is the first step towards developing an efficient implementation of AR$3$.   We refer the reader to Chapters 4 and 12 in \cite{cartis2022evaluation} for further discussion of the challenges and benefits of using high-order methods.}

\subsection{Minimizing the Quartically-regularised Cubic Polynomial $m_3$}
\label{section_1.1}
In the case $p = 3$, the sub-problem $m_3$ is a possibly nonconvex, quartically-regularised multivariate polynomial,
\begin{equation}
\tag{AR$3$ Model}
 m_{3}(x_k, s) = f(x_k) + \nabla_x f(x_k)^Ts + \frac{1}{2} s^T\nabla_x^2 f(x_k) s + \frac{1}{6}  \nabla_x^3 f(x_k) [s]^3+ \frac{ \sigma_k}{4} \|s\|^4.
  \label{ar3 model}
\end{equation} 
While Nesterov proved that globally minimizing a cubic polynomial over a ball is NP-hard \cite{burachik2021steklov, luo2010semidefinite, nesterov2003random}, we aim to efficiently locate a local minimum of $m_3(s)$, which according to \cite{cartis2020concise}, is sufficient to ensure the good complexity of AR$3$.
For notational simplicity, we fix $x_k$ and write $m_3(x_k, s)$ as
\begin{equation}
 m_3(s) =f_0 + g^Ts+ \frac{1}{2} s^T H s + \frac{1}{6}  T [s]^3+ \frac{\sigma}{4}  \|s\|^4,
 \label{m3}
\end{equation}
where $f_0 = f(x_k) =  m_3(x_k, 0) \in \R$, $g = \nabla_x f(x_k) \in \R^n$, $H = \nabla^2_x f(x_k) \in \R^{ n \times n}$ and $T = \nabla^3_x f(x_k) \in \R^{ n \times n \times n}$. 
We also denote the fourth-order Taylor expansion of $m_{3} (s^{(i)}+s)$ at $s^{(i)}$ as
\begin{eqnarray*}
M(s^{(i)}, s)  := f_i + g_i^T s + \frac{1}{2}s^T H_i s +\frac{1}{6}T_i[s]^3 + \frac{\sigma}{4}\|s\|^4,
\end{eqnarray*} 
where $f_i =  m_3(s^{(i)}) \in \R$, $g_i = \nabla m_3(s^{(i)}) \in \R^n$, $H_i = \nabla^2 m_3(s^{(i)}) \in \R^{ n \times n}$ and $T_i = \nabla^3 m_3(s^{(i)}) \in \R^{ n \times n \times n}$ and $\nabla$ denotes the derivative with respect to $s$.  
For all $s^{(i)}$,  the Hessian $H_i$ is a symmetric matrix, and $T_i \in \R^{n\times n \times n}$ is a symmetric third-order tensor\footnote{A symmetric tensor means that the entries are invariant under any permutation of its indices.}. The first- and second-order derivatives of $M(s^{(i)}, s)$ are given by $\nabla M(s^{(i)}, s) = g_i + H_i s +\frac{1}{2}T_i[s]^2 + \sigma\|s\|^2 s$, and $\nabla^2 M(s^{(i)}, s) = H_i +T_i[s] + \sigma (\|s\|^2 I_n + 2 ss^T)$, respectively. Here, $I_n \in \R^{n \times n}$ denotes the identity matrix, and the matrix $ss^T$ is a rank-one matrix. We refer to the model $M(s^{(i)}, s)$ as locally convex at $s^{(i)}$ if $H_i \succeq 0$, and it is locally nonconvex at $s^{(i)}$ if the smallest eigenvalue of $H_i$, denoted as $\lambda_{\min}[H_i]$, is negative.

Since $m_3(s)$ is a fourth-degree multivariate polynomial, the fourth-order Taylor expansion is exact. Therefore, we have
\begin{eqnarray*}
M(s^{(i)}, s) =m_3(s+s^{(i)}),   \qquad  \text{ and } \qquad \min_{s\in \R^n}m_3(s) = \min_{s\in \R^n}M(s^{(i)}, s). 
\end{eqnarray*}
Our aim is to find an iterate $s^{(i)}$ such that 
\begin{equation}
\label{opt cond 2}
    \big\|g_i\big\| \leq \epsilon,
\end{equation}
where $\epsilon \in (0, 1)$ is an arbitrary user-set tolerance. {This implies that $\liminf_{i \rightarrow \infty}g_i = 0$, and at least one limit point of the sequence of iterates $\{s^{(i)}\}_{i \ge 0}$ is a first-order stationary point of $m_3$.}

There are some existing algorithmic and numerical proposals involving related AR$3$ sub-problems. Birgin et al \cite{birgin2017use} present a variant of the AR$3$ algorithm and its implementation, comparing its performance with that of AR$2$. Schnabel et al \cite{chow1989derivative, schnabel1971tensor, schnabel1991tensor} also consider solving unconstrained optimization using third-order tensor local models.

Nesterov proposed several methods for minimizing a convex $m_3$ model, including second-order methods utilizing Bregman gradient methods under relative smoothness properties \cite{Nesterov2006cubic, Nesterov2021implementable}, difference approximation of the third-order tensor term \cite{Nesterov2021superfast}, and high-order proximal-point operators \cite{Nesterov2020inexact, Nesterov2021inexact}.
More recently, Nesterov gave a linearly convergent second-order method for minimizing convex quartically-regularised cubic polynomials \cite{Nesterov2022quartic}. We refer to this method as Nesterov's method. At each iteration, Nesterov approximates the third-order tensor term by a linear combination of second-order derivatives and fourth-order regularization. Then, the algorithm minimizes $m_3$ by minimizing a sequence of convex quadratic models with quartic regularization, iteratively. To the best of our knowledge, there are no solvers specifically designed for the (local) minimization of nonconvex quartically regularized cubic polynomials.

The key idea presented in this paper involves the generalization of Nesterov's method to the minimization of nonconvex $m_3$ {models}. We refer to this framework as the \textbf{Quadratic Quartic Regularization Method} (QQR method). In summary, the QQR method minimizes convex or nonconvex AR$3$ sub-problems of the form \eqref{m3} by minimizing a series of quadratic polynomials (that maybe nonconvex) with quartic regularization. 

\textbf{The paper makes three key contributions:}
\begin{itemize} \itemsep0em 
    \vspace{-0.2cm}\item We proposed an algorithmic framework, the QQR method, for efficiently minimizing possibly, nonconvex, quartically-regularised cubic polynomials. To the best of our knowledge, this is the first solver designed to minimize the AR$3$ sub-problem in \eqref{m3} for a nonconvex objective function. QQR method can be seen as a nonconvex generalization of Nesterov's method proposed in \cite{Nesterov2022quartic}, which was originally designed for convex $m_3$ functions. When applied to convex $m_3$ functions, the QQR method exhibits a similar linearly convergent behaviour as Nesterov's method. 
        
    \item Theoretically, we show that the QQR method achieves better complexity bounds than the cubic regularization method in {locally strictly convex and locally nonconvex} cases. For locally convex iterations, the QQR method either gives a linear rate of convergence for the norm of the gradient or an error reduction of at least $\mathcal{O}(\epsilon^{4/3})$, where $\epsilon$ represents the accuracy of the first-order optimality condition \eqref{opt cond 2}. In locally nonconvex iterations, the QQR method achieves an error reduction of at least $\mathcal{O}(\epsilon)$. 

    \vspace{-0.2cm}\item We propose two practical variants of the QQR method. In our numerical experiments, QQR performs competitively with state-of-the-art approaches such as ARC (also known as AR$p$ with $p=2$ \cite{cartis2011adaptive}) {for the minimization of $m_3$}. It achieves better results by finding a minimum with fewer algorithmic iterations or function evaluations, or by finding a lower minimum. 
\end{itemize}

\subsection{Related Work Regarding Minimization of Quartic Polynomials}

In addition to the existing works already mentioned 
in Section \ref{section_1.1},
in some special cases, such as when the polynomial is convex or can be expressed as a sum of squares (SOS) of polynomials~\cite{lasserre2001global, lasserre2015introduction, laurent2009sums}), semidefinite programming (SDP) methods (such as the Lasserre hierarchy) can be shown to converge to a global minimizer of the quartic polynomial. However, alternative methods are required when polynomials cannot be expressed in SOS form.
One approach proposed by Qi \cite{qi2004global} focuses on normal quartic polynomials with a positive definite fourth-degree coefficient tensor. The authors develop a global descent algorithm to locate the global minimum of normal quartic polynomials in the case where $n=2$. Another strategy is to utilize branch-and-bound algorithms, which involve recursively partitioning the feasible region and constructing nonconvex quadratic or cubic lower bounds. These algorithms, such as those presented in \cite{cartis2015branching, fowkes2013branch, horst2013global, kvasov2009univariate, neumaier2004complete}, aim to identify the global minimum through bounding procedures. In recent work by Burachik et al. \cite{arikan2020steklov, burachik2021steklov}, a trajectory-type method is proposed for global optimization of multivariate quartic normal polynomials. This approach convexifies the quartic polynomial and generates a path originating from the minimizer of the convexified function. While convergence is proven for univariate quartic normal polynomials, the results for multivariate cases are only supported by numerical experiments. Despite the advancements in these methods, they all face challenges associated with the curse of dimensionality. The computational complexity can grow exponentially with the polynomial's number of variables, rendering certain methods impractical for large-scale problems. Local first and second-order methods devised for locally minimizing $f(x)$ can directly be applied to minimizing $m_3$; in particular, AR$2$ is a suitable choice. We will discuss this and other similar methods in more detail later in the paper.
During the preparation of this revision, there has already been some follow-up research. One such study combines the polynomial sum-of-squares framework with adaptive regularization techniques for tensor methods, providing the first tensor methods with provably polynomial work per iteration, in the convex case \cite{ahmadi2024higher}.  Subsequently, \cite{cartis2024global} extended this approach to the nonconvex case and  analyzed its non-asymptotic rate of convergence and complexity  for finding approximate stationary points. Another contribution  \cite{zhu2023cubic} focuses on  minimizing \eqref{m3} by employing a sequence of local quadratic models that incorporate simple cubic and quartic terms.  Still, the current approach here stands on its own, as having  efficient  subproblem reformulations.

\subsection{Motivation for Our Work}

Our QQR method is inspired by the works of Cartis, Gould and Toint~\cite{cartis2022evaluation} and Nesterov~\cite{Nesterov2022quartic}. In this section, we motivate using the quadratic quartically-regularised model and discuss optimality conditions for this model.

\textbf{The choice of a quadratic quartically-regularised model} The first question to address regarding the QQR method is why we choose a quadratic quartically-regularised model to minimize the AR$3$ sub-problem \eqref{m3}. Consider the second-order model with an $r$th-power regularization:
\begin{equation*}
m_{2}^r(s) = \hat{f_0} + \hat{g}^Ts + \frac{1}{2}  {s^T\hat{H}s} + \frac{1}{r} \sigma \|s\|^{r}, 
\end{equation*}
where $r > 2$ and $\hat{f_0} \in \R$, $\hat{g} \in \R^n$ and $\hat{H} \in \R^{n \times n}$ are scalar, vector and symmetric matrix coefficients. $m^3_2(s)$ has the same form as the ARC sub-problem and $m^4_2(s)$ has the same form as our quadratic quartic regularization model. The minimization of the potentially nonconvex $m_{2}^r(s)$ has been widely researched as part of the (adaptive) cubic regularization framework~\cite{cartis2011adaptive, dussault2018arcq, martinez2017cubic, kohler2017sub, Nesterov2006cubic}. 
Theorem \ref{thm: arc} summarizes the necessary and sufficient optimality conditions for the global minimizer of $m^r_2(s)$. 
\begin{theorem}
 (Theorem 8.2.8~\cite{cartis2022evaluation}) 
Let $r \ge 3$, any global minimizer of $m_{ 2}^r(s)$, $s_*^r$, satisfy
$
    (\hat{H} + \hat{\lambda}_* I_n) s_*^r = -\hat{g},
$
where $I_n \in \R^{n \times n}$ is the identity matrix, $\hat{\lambda}_* \ge 0$ , $\hat{H} + \hat{\lambda}_* I_n \succeq 0 $, and 
$\hat{\lambda}_* = \sigma \|s_*^r\|^{r-2}. $ If $\hat{H} + \hat{\lambda}_* I_n$ is positive definite, then $s_*^r$ is unique.
\label{thm: arc}
\end{theorem}

Theorem \ref{thm: arc} converts the task of identifying the global minimum of $m_2^r$ into a problem involving the solution of a nonlinear equation combined with a matrix linear system. 
 In broad terms, we use the matrix system $(\hat{H} + \hat{\lambda}_* I_n) \hat{s}_* = -\hat{g}$ to express $s^*$ as a function of $\hat{\lambda}_*$, such that $\hat{s}_* = \hat{s}(\hat{\lambda}_*)$. We obtain an approximate global minimizer $(\hat{s}_*, \hat{\lambda}_*)$ by solving the nonlinear equation $\hat{\lambda}_* = \sigma \|\hat{s}(\hat{\lambda}_*)\|^{r-2}$.  There are many scalable methods for minimizing the sub-problems in $m_2^r$~\cite{cartis2011adaptive, martinez1994local, lucidi1998some, hsia2017theory}. Specifically, \cite{cartis2011evaluation, cartis2022evaluation} describe a factorization-based minimization solver for quadratic problems with fourth-order regularization, which we refer to here as the \texttt{MS} algorithm (due to its connection to Mor\'e-Sorensen techniques as described next). The \texttt{MS} algorithm finds (an approximate)  global minimizer of both convex and nonconvex $m_2^r$, by applying Newton’s method for root finding to a secular equation together with a Cholesky factorization of the optimality conditions that hold at a global minimizer of $m_2^r$; this  is in the same vein as the Mor\'e-Sorensen approach for solving the trust-region subproblem. More details can be found in  Ch 8~\cite{cartis2022evaluation}.

\textbf{Nesterov's method for convex $m_3$:} A summary of Nesterov's results is given in Remark \ref{Nesterov result}. 

\begin{tcolorbox}[breakable, enhanced]
\begin{remark}
\label{Nesterov result} Assume that $m_3(s)$ is \textbf{strictly convex} for all $s \in \R^n$. Then, 
\begin{itemize} \itemsep0em 
\vspace{-0.2cm}\item \textbf{Convex Upper and Lower Bounds} \cite[Corollary 1]{Nesterov2022quartic}:  For any $\tau \in (\frac{1}{3}, \frac{1}{2})$, let $ \alpha^{\pm} = [\alpha_1^{\pm}, \alpha_2^{\pm}] = [1 \pm \frac{1}{3\tau}, 1\pm 2\tau]$, such that
\begin{equation}
M^{\mathcal{N}}_{\alpha^{\pm}}(s^{(i)},s) = f_i + g_i^Ts+  \frac{1}{2} (1 \pm \frac{1}{3\tau}) s^TH_is  + \frac{\sigma}{4}( 1\pm 2\tau)\|s\|^4
\label{nesterov upper lower}
\end{equation}
Then, both $M^{\mathcal{N}}_{\alpha^{\pm}}(s^{(i)},s)$ are convex and $
M^{\mathcal{N}}_{\alpha^-}(s^{(i)},s) \le M(s^{(i)},s) \le M^{\mathcal{N}}_{\alpha^+}(s^{(i)},s) $
 for all $s \in \R^n$, and the equality only occurs when $s \equiv 0 $.
    \vspace{-0.2cm}\item \textbf{Linear Convergence Result} \cite[Theorem 3]{Nesterov2022quartic}: Fix $\tau =\tau_{\mathcal{N}} := \frac{1}{6}\sqrt{(3+\sqrt{33})}$. We update the iterates $s^{(i)}$ by 
    $s_+^{(i)} = \argmin_{s\in \R^n} M_{\alpha^+}(s^{(i)}, s),$ and $s^{(i+1)} = s^{(i)} + s_+^{(i)}.$
   Then, $m_3(s^{(i)})$ satisfies
    \begin{equation*}
        m_3(s^{(i)}) - m_{*}\le \left(1-\mu_{\mathcal{N}} \right)^{i} ( m_3(s^{(0)}) - m_{*})  
    \end{equation*}
    where $m_* = \min_{s\in\R^n} m_3(s)$ and $\mu_{\mathcal{N}} = 0.193$.
\vspace{-0.2cm}\item $M^{\mathcal{N}}_{\alpha^+}$ ensures the monotonic model decrease in $m_3(s^{(i)})$ as $i$ increases.  $M^{\mathcal{N}}_{\alpha^-}$ is used in the complexity proof to construct $M^{\mathcal{N}}_\mu$.  The upper and lower bounds are linked by $M^{\mathcal{N}}_\mu$ defined as
\begin{equation}
M^{\mathcal{N}}_\mu (s^{(i)}, s) := f_i + g_i^T s + \mu_{\mathcal{N}}^{-1} \bigg( \frac{1}{2} - \frac{1}{6\tau_{\mathcal{N}}}
\bigg)s^TH_is + \mu_{\mathcal{N}}^{-3} \bigg(1-2\tau_{\mathcal{N}} \bigg)\frac{\sigma}{4}\|s\|^4.
\label{nesterov mu bound}
\end{equation}
\end{itemize}
\end{remark}
\end{tcolorbox}
Based on the above approach,  we make three key remarks.  
\begin{itemize} \itemsep 0em 

    \vspace{-0.2cm}\item In \cite{Nesterov2022quartic}, $  M^{\mathcal{N}}_{\alpha^\pm}$ is a family of $\tau$-parameterized upper and lower bounds. For convex quartically-regularised cubic polynomials minimization,  user-chosen $(\tau_\mathcal{N}, \mu_\mathcal{N})$ are sufficient to achieve the linear convergence result. However, in the nonconvex case, it is necessary to incorporate some adaptive algorithm parameters to ensure the convergence of $m_3(s^{(i)})$.
    \vspace{-0.2cm}\item Although the convergence proof for Nesterov's algorithm only applies to globally convex quartic polynomials, we extend it to handle quartic polynomials with only local convexity. We show that under suitable local convexity assumptions, both the linear rate of convergence and upper/lower bound results for convex $m_3$ can be maintained.
\end{itemize}

 These observations inspired us to create a QQR framework that minimizes generally nonconvex quartically-regularised cubic polynomials.  We demonstrate that, for convex $m_3$ problems, QQR recovers both of Nesterov's results, namely the upper/lower bound results and the linear convergence result.\footnote{It is worth noting that Nesterov's method can be applied directly to the minimization of a convex objective function (not just $m_3$). Similarly, if the objective function satisfies suitable Lipschitz conditions and is bounded below, the QQR method can also be used for local minimization of a potentially nonconvex objective function. However, our interest here is not to create a general new second-order method, but to specifically attempt to efficiently minimize $m_3$.} Example \ref{QQR n=1} is a simple illustration of the QQR algorithm for a nonconvex univariant quartic polynomial. 

\begin{example}
QQR gives strict descent in $m_3(s^{(i)})$ in each iteration. In locally nonconvex iterations (i.e., $\lambda_{\min}[H_i]<0$), we also have a nonconvex upper and a nonconvex lower bound. 
\begin{figure}[ht!]
  \begin{center}
    \includegraphics[width=12cm]{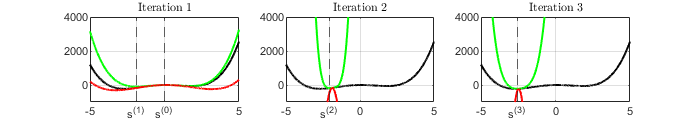}
    \caption{\small  Minimizing $m_3(s) = 10s-50s^2+5s^3+5s^4$. $m_3$, $M_{\alpha^+}$, $M_{\alpha^-}$ plotted in black, green, and red respectively. $M_{\alpha^-}$ is bounded below on each iteration but the lower bound is outside the range of the plots. }
    \label{QQR n=1}
  \end{center}
\end{figure}
    \label{example 1}
\end{example}

The paper is structured as follows. Section \ref{section_2} serves as the core of the paper, introducing the QQR algorithmic framework and two main variants. In this section, we explore the connection between the QQR methods and Nesterov's method, specifically for convex $m_3$ or locally convex iterations. We provide a comprehensive theoretical analysis of the complexity and convergence properties of the QQR framework. In Section \ref{sec tech proof}, we investigate the technical lemmas required for the convergence of the QQR algorithm. These theoretical findings form the basis for the two implementations of the QQR method presented in Section \ref{sec_implementation}.  We provide numerical examples in Section \ref{sec_numerics}, comparing our proposed QQR variants with the state-of-the-art cubic regularization method. Finally, Section \ref{Conclusion} offers concluding remarks.

\section{The QQR Method}
\label{section_2}

QQR minimizes the AR3 sub-problem (i.e., $m_3$ as defined in \eqref{m3}) by globally optimizing a series of potentially nonconvex quadratic models with quartic regularization. 
The QQR algorithm divides iterations into three distinct cases: \textbf{locally strictly convex}, \textbf{locally nearly convex}, and \textbf{locally nonconvex}, {which we refer to as Cases 1, 2, and 3, respectively}.  Depending on the local convexity of $m_3$ at each iteration $i$, the parameters for the quadratic term and quartically regularized term are specifically chosen. Algorithm \ref{TQQR} gives the framework of the QQR algorithm. 

\begin{algorithm}[!ht]
\caption{\small General QQR Framework}
\label{algo1}
\small{
\textbf{Initialization}: Initial guess at $s^{(0)} = 0 \in \R^n$, $i:=0$. Set the tolerance for the first-order optimality $\epsilon$. 
\\ Set $\lambda_{c} > 0 $; $\lambda_{c}$ is a user-defined constant used to distinguish among the three distinct cases: locally nonconvex, locally strictly convex, and locally nearly convex. Usually, set $\lambda_c:= \epsilon^{1/3}$ and {$\eta \in (0, \frac{1}{3})$}.

\While{$\|g_i\|  >\epsilon $}
{\textbf{Step 1: QQR Model Selection.}
Set  $
M_{\alpha^+}(s^{(i)}, s) : = f_i + g_i^T s + \frac{\alpha_1^+}{2} s^T \big(H_i  + \tilde{p} I_n\big) s+ \frac{\sigma}{4}\alpha_2^+ \|s\|^4 
$
where $\alpha_1^+>0$, $\alpha_2^+ > 0$ and $\tilde{p} \ge 0$ are chosen specifically such that 
\begin{eqnarray*}
 M_{\alpha^+}(s^{(i)}, s) \text{ is a (local) upper bound/approximation for }  m_3(s) \text{ at } s = s^{(i)}. 
\end{eqnarray*}

The method for selecting $\alpha^+ = {(\alpha_1^+, \alpha_2^+)}$ will be discussed in later sections. Here, we provide a summary of the key features of these constants. 
\begin{enumerate}
    \vspace{-0.1cm}\item\textit{Locally nonconvex iteration:} i.e.,  $\lambda_{\min}[H_i] \le -\lambda_{c}$. Choose $0< \alpha_1^+ \le 1$  and $\tilde{p} =0$. 
    \vspace{-0.2cm}\item \textit{Locally strictly convex iteration:}  i.e., $\lambda_{\min}[H_i] \ge \lambda_c$. 
    Choose $\alpha_1^+ > 1$ and $\tilde{p} =0$. 
    \vspace{-0.2cm}\item \textit{Locally nearly convex iteration:} i.e.,  $ \big| \lambda_{\min}[H_i] \big| < \lambda_c$. Choose $\alpha_1^+=1$ and $  \tilde{p} > 0$. 
\end{enumerate}
\textbf{Step 2: Iterative Minimization.}
\begin{equation}
    \qquad  s_+^{(i)} = \argmin_{s\in \R^n} M_{\alpha^+}(s^{(i)}, s), \qquad  s^{(i+1)} = s^{(i)} + s_+^{(i)} , \qquad i: = i+1.
    \label{iterates}
\end{equation}
\\ 
\textbf{Step 3: Test for Acceptance.} Let
\begin{equation}
\tag{Ratio Test}
\rho(s_+^{(i)}) := \frac{m_3(s^{(i)}) - m_3(s^{(i)}+s_+^{(i)})}{m_3(s^{(i)}) -  M_{\alpha^+}(s^{(i)}, s_+^{(i)})}.
\label{ratio test}
\end{equation}
If $\rho(s_+^{(i)}) \ge \eta$,  set $s^{(i+1)} := s^{(i)}+{s_+^{(i)}}$. Otherwise, set $s^{(i+1)} := s^{(i)}$ and adjust $\alpha_1^+, \alpha_2^+$ and $\tilde{p}$ adaptively. Let $i=i+1$.
}
}
\label{TQQR}
\end{algorithm}

The section is organized as follows. In Section \ref{sec: locally strictly convex Iterations}, we present QQR variant 1 with a single adaptive model parameter and discuss the relationship between the QQR methods and Nesterov's method. In Section \ref{sec: theory QQR}, we present QQR variant 2 with two adaptive model parameters. In Section \ref{Proof of Complexity}, we provide a case-by-case theoretical analysis of the convergence and complexity of QQR methods. We show that QQR methods achieve better complexity than the cubic regularization method during locally strictly convex or nonconvex iterations.

\subsection{QQR Variant 1: Single Adaptive Model Parameter}
\label{sec: locally strictly convex Iterations}

Similarly to Nesterov's approach, we construct a QQR variant that adaptively adjusts one parameter. This parameter governs the interplay between the quadratic term and the quartic term.

\begin{tcolorbox}[breakable, enhanced]

At the $i$th iteration of the QQR Framework (Algorithm \ref{algo1}), we iteratively minimize the QQR variant 1 model (\texttt{QQR-v1}) given by:
\begin{eqnarray}
M_\mu (s^{(i)}, s) : = f_i + g_i^T s + \frac{\alpha_1^\mu}{2} s^T\big[H_i + \tilde{p}I_n \big]s + \underbrace{\mu^{-\kappa} \bigg(1-2(1-q)\tau\bigg)}_{:=\alpha_2^\mu}\frac{\sigma}{4}\|s\|^4
\label{qqr v1}
\end{eqnarray}
where
\begin{eqnarray*}
\text{Locally strictly convex iteration:} &&\alpha_1^\mu  = \mu^{-\iota} \bigg( 1- \frac{1}{3\tau}\bigg),  \qquad \tilde{p} = 0 
\\ \text{Locally nearly convex iteration:} &&\alpha_1^\mu  = \mu^{-\iota}  \bigg( 1- \frac{1}{3\tau}\bigg) ,  \qquad \tilde{p} > 0 
\\ \text{Locally nonconvex iteration:} &&\alpha_1^\mu = 1,  \qquad  \tilde{p} = 0. 
\end{eqnarray*}
$\tau, q, \kappa, \iota$ and $\tilde{p}$ are  fixed constants\footnote{Details regarding choosing these parameters are given in Section \ref{sec: algo for QQRv1}.} and (only) $\mu$ is adaptively adjusted\footnote{For theoretically suitable choices of $\mu$, see Remark \ref{theoretical requirement for mu} and Corollary \ref{link to qqr1}; for practical $\mu$, see Section \ref{sec: algo for QQRv2}.}. 
\end{tcolorbox}

\begin{remark}
$\tau, \iota, \kappa>0$ and $q \in [0, 3)$ are initialized as fixed constants.
For a sufficiently small positive $\mu$, we can guarantee  $\alpha_1^\mu > 1$ for locally strictly convex iterations (see Algorithm \ref{TQQR}). 
In practice, we achieve the desired small value of $\mu$ by iteratively reducing its magnitude.

\label{qqr v1 parameter}
\end{remark}

\begin{remark}
For a strictly convex $m_3$, if we set $\tau = \tau_\mathcal{N}$, $q=0$, $\tilde{p}=0$, $\iota = 1$, $\kappa =3$, then the QQR-v1 model in \eqref{qqr v1} coincides with the model proposed by Nesterov in \eqref{nesterov mu bound}; namely,  $M_\mu = M^{\mathcal{N}}_\mu$. However, it is important to note that when minimizing a convex $m_3$, one can choose fixed values of $(\tau, \mu)$ at $(\tau_\mathcal{N}, \mu_\mathcal{N})$. Therefore, in Nesterov's method, only the upper bound $M^{\mathcal{N}}_{\alpha^+}$ needs to be iteratively minimized. In contrast, since we are minimizing a nonconvex $m_3$, we need to adaptively adjust $\mu$, and thus we iteratively minimize $M_\mu$ to allow for changes in the model based on convexity.
\label{remark relationship of qqr1 and nesterov}
\end{remark}





In the following sections, we show that in locally strictly convex iterations when $\lambda_{\min}[H_i] > \lambda_c$, \texttt{QQR-v1} locally recovers the upper/lower bounds in Section \ref{sec: recover result 1}. Additionally, we establish that \texttt{QQR-v1} achieves a linear rate of convergence within a region (Section \ref{sec: recover result 2}).

\subsubsection{Upper and Lower Bounds in $\mathcal{D}^{(i)}(q)$}
\label{sec: recover result 1}

We first define the region within which \texttt{QQR-v1} achieves a linear rate of convergence. 
{\begin{definition}
 Let $i$ be a locally strictly convex iteration of \texttt{QQR-v1} so that $ \lambda_{\min}[H_i]  \ge \lambda_c >0$. 
 For a(ny) $q \ge 0$, define a region $\mathcal{D}^{(i)}(q) \subseteq \R^n$ around $s^{(i)}$, such that
 for $s \in \mathcal{D}^{(i)}(q)$, the following conditions are satisfied: 
\begin{enumerate}
    \item $ s^T \nabla^2 m_3(s^{(i)}+ s) s \ge q \|s\|^4$ and $ s^T \nabla^2 m_3(s^{(i)}- s) s \ge q \|s\|^4$;
    \item For any $s\in \mathcal{D}^{(i)}(q)$, we also have $\tau s\in \mathcal{D}^{(i)}(q)$.
\end{enumerate}
\label{def dq}
\end{definition}

Clearly, for any $q \ge 0$, $s = 0$ satisfies both conditions for $\mathcal{D}^{(i)}(q)$, therefore $\mathcal{D}^{(i)}(q)$ is well-defined and nonempty. Also, the first condition in Definition \ref{def dq} indicates that, if $s \in \mathcal{D}^{(i)}(q)$, then $-s \in \mathcal{D}^{(i)}(q)$. 

\begin{lemma} There exists  $\delta_i>0$
such that $\bar{\mathcal{B}}(0, \delta_i) := \{s \in \R^n, \|s\| \le \delta_i\}  \subseteq \mathcal{D}^{(i)}(q)$. 
\end{lemma}
\begin{proof}
Since $H_i = \nabla^2 m_3(s^{(i)}) \succeq \lambda_c I_n$, there exists   $\hat{\delta_i}>0$  such that for all $s \in \bar{\mathcal{B}}(0,\hat{\delta_i})$, we have $ \nabla^2 m_3(s^{(i)}+s) \succeq \lambda_c I_n$. If $q=0$, set $\delta_i := \hat{\delta}_i $; otherwise choose $\delta_i := \min\{\frac{\lambda_c}{q},  \hat{\delta}_i\}$. For all $s \in \bar{\mathcal{B}}(0, \delta_i)$, we have
 $s^T  \nabla^2 
 m_3(s^{(i)}+s) s  \ge \lambda_c \|s\|^2 \ge q \|s\|^4$. 
Clearly, the latter inequalities also hold  for 
$-s$ and $\tau s$ for $\tau \in (0, 1)$. Therefore, $s \in \mathcal{D}^{(i)}(q)$ and $ \bar{\mathcal{B}}(0, \delta_i) \subseteq \mathcal{D}^{(i)}(q)$.
\end{proof}

\begin{remark} ($\mathcal{D}^{(i)}(q)$ and Convexity) If $m_3$ is strictly convex for all $s \in \R^n$, then $s^T \nabla^2 m_3(s^{(i)} + s) s > 0$ for all $s, s^{(i)} \in \R^n$ implying that $\mathcal{D}^{(i)}(0) = \mathbb{R}^n$. More generally, the size of the region $\mathcal{D}^{(i)}(q)$ depends on the value of $q$. When we set $q = 0$, condition 1 in Definition \ref{def dq} is less stringent, resulting in a larger region $\mathcal{D}^{(i)}(q)$ around $s^{(i)}$. Conversely,  larger values of $q$ lead to smaller regions $\mathcal{D}^{(i)}(q)$. 
\end{remark}
}


 Within the region $\mathcal{D}^{(i)}(q)$, we can provide a convex upper and lower bound for $m_3$,  Theorem \ref{thm: Ada Nesterov convex} shows.

\begin{theorem} (Upper and Lower Bounds in $\mathcal{D}^{(i)}(q)$) Let $i$ be a locally strictly convex iteration of \texttt{QQR-v1} such that $ \lambda_{\min}[H_i] \ge \lambda_c >0$. {Fix  $q \in [0, 3)$ and $\tau$  such that
\begin{eqnarray}
\frac{1}{3} < \tau < \min \bigg\{\frac{1}{2(1-q/3)}, 1 \bigg\}.
\label{q and tau}
\end{eqnarray}}Let
\begin{eqnarray}
    M_{\alpha^{\pm}}(s^{(i)}, s) = f_i + g_i^T s +   \frac{1}{2}\bigg(1\pm\frac{1}{3\tau}\bigg) s^TH_i s  +\bigg[1\pm 2\big(1-\frac{q}{3}\big)\tau\bigg] \frac{\sigma}{4}\|s\|^4
    \label{M alpha tau}
\end{eqnarray}
Then, 
\begin{description}
    \vspace{-0.2cm}\item 1)  $M_{\alpha^{\pm}}(s^{(i)}, s)$ are convex functions for all $s\in \R^n$. 
    \vspace{-0.2cm}\item 2) $M_{\alpha^-}(s^{(i)},s) \le M(s^{(i)},s) \le M_{\alpha^+}(s^{(i)},s) $ for all $s \in \mathcal{D}^{(i)}(q)$, with equality if and only if $s = 0 $.
\end{description}
\label{thm: Ada Nesterov convex}
\end{theorem}

\begin{proof}
\textbf{Result 1): } For $\frac{1}{3} < \tau < \min\{\frac{1}{2(1-q/3)}, 1\}$, we have $\frac{1}{2}\pm\frac{1}{6\tau} > 0$ and $ 1\pm 2(1-q/3)\tau> 0$.
Thus, \small{$\nabla^2 M_{\alpha^\pm}(s^{(i)},s) =   \underbrace{\bigg(\frac{1}{2}\pm\frac{1}{6\tau}\bigg) H_i }_{\succ 0}  + \underbrace{\frac{\sigma}{4} \bigg[ 1\pm 2\big(1-\frac{q}{3}\big)\tau \bigg] \bigg(\|s\|^2 I_n +2ss^T\bigg)}_{\succ 0}\succ 0. $} \normalsize
    
\textbf{Result 2): } 
For any $s\in \mathcal{D}^{(i)}(q)$, we have $\pm \tau s\in \mathcal{D}^{(i)}(q)$. Specifically, 
\begin{eqnarray}
(\tau s)^T \nabla^2 m_3(s^{(i)} \pm  \tau s) (\tau s) - q \| \tau s\|^4  >0. 
\label{mid step 1}
\end{eqnarray}
An (exact) second-order Taylor expansion of $\nabla^2 m_3(s^{(i)} \pm \tau s)$ implies
\begin{eqnarray}
\nabla^2 m_3(s^{(i)} \pm  \tau s) = H_i \pm  \tau T_i[ s] + 3\sigma  \tau^2 \| s\|^2I_n,  
\label{mid step 5}
\end{eqnarray}
Substituting \eqref{mid step 5} into \eqref{mid step 1}, we deduce that
$s^T \big[ H_i  \pm  \tau T_i[s] + (3-q) \sigma  \tau^2 \| s\|^2 \big]s > 0$ for $0 < \tau <1$. Now, consider any  $s\in \mathcal{D}^{(i)}(q)$, 
\small{
\begin{eqnarray*} 
M_{\alpha^+}(s^{(i)},s) - M(s^{(i)},s)   = \frac{1}{6\tau} s^TH_i s - \frac{1}{6}T_i[ s]^3 +  \frac{3-q}{6}\sigma\tau \| s\|^4 = \frac{1}{6\tau} s^T\bigg[  H_i- \tau T_i[s] +  (3-q) \sigma \tau^2 \| s\|^2  \bigg]s 
>0.
\end{eqnarray*}}
\normalsize
 The proof for $M_{\alpha^-}(s^{(i)},s)$ follows similarly.
\end{proof}

\begin{remark}
When $q = 0$, the upper and lower bounds $M_{\alpha^{\pm}}$ in \eqref{M alpha tau} coincides with $M^{\mathcal{N}}_{\alpha^{\pm}}$ in \eqref{nesterov upper lower} in Nesterov's model. {Additionally, choosing $0 \le q \le 3$ ensures that the regularization parameter of the quartic  term stays positive, i.e. $1\pm 2\big(1-\frac{q}{3}) \ge 0$.}
\end{remark}

\subsubsection{Linear Convergence of \texttt{QQR-v1}  in $\mathcal{D}^{(i)}(q)$}
\label{sec: recover result 2}

{In this subsection, we prove that \texttt{QQR-v1} converges linearly to the minimum of the quartically-regularised cubic polynomial in $\mathcal{D}^{(i)}(q)$ over locally strictly convex iterations.} To prove this, we interpret $M_\mu(s^{(i)}, s)$ in \eqref{qqr v1}  as an adjustable bound that increases from $M_{\alpha^-}$ to $M_{\alpha^+}$ as $\mu$ decreases from $1$ to $0$. It is clear that when $\mu = 1$, $M_\mu$ coincides with our local lower bound $M_{\alpha^-}$, and as $\mu \rightarrow 0$, $M_\mu(s^{(i)}, s)$ tends to infinity for all $s \in \R^n /\{0\}$. Before presenting the complexity proof for a locally strictly convex iteration, we provide a technical lemma for $M_\mu(s^{(i)}, s)$.

\begin{lemma}
Let $M_\mu(s^{(i)}, s)$ be defined  in \eqref{qqr v1} with $\iota > 0$, $\kappa > 0$, $\tilde{p} = 0$ and $(q, \tau)$ satisfying \eqref{q and tau}. If $\mu$ satisfies 
\begin{eqnarray}
    \mu \le \min \bigg\{\bigg[\frac{3\tau - 1}{3\tau + 1}\bigg]^{1/\iota}, \bigg[\frac{1 -2(1-q/3)\tau}{1 +2(1-q/3)\tau} \bigg]^{1/\kappa} \bigg\},
    \label{cond on mu}
\end{eqnarray}
then $M_\mu(s^{(i)}, s) \ge M_{\alpha^+}(s^{(i)}, s)$
for all $s \in \R^n$. 
\label{cond alpha convex} 
\end{lemma}

\begin{proof}  
The condition on $\mu$ ensures $  \mu^{-\iota}(\frac{1}{2} - \frac{1}{6\tau}) \ge \frac{1}{2} + \frac{1}{6\tau} \ge 0$ and $  \mu^{-\kappa}[1-2(1-q/3)\tau] \ge [1+2(1-q/3)\tau] \ge 0$. Additionally, using the fact that $s^TH_i s > 0$ for any $s \in \mathbb{R}^n$ yields the desired result.
\end{proof}

Now, we present the theorem that establishes the linear reduction in the function value in a locally strictly convex iteration.

\begin{theorem} 
\label{thm: convex complexity result}
(Linear Rate of Convergence in $\mathcal{D}^{(i)}(q)$.) 
Assume that $i$ is a locally strictly convex iteration of \texttt{QQR-v1}  satisfying $\lambda_{\min}[H_i] \ge \lambda_c$ and $\|g_i\| > \epsilon$, and construct $M_{\alpha^+}(s^{(i)}, s)$ as  in Theorem \ref{thm: Ada Nesterov convex}. We update the iterates via iterative minimization in $\mathcal{D}^{(i)}(q)$, 
\begin{equation}
 s_q^{(i)} = \argmin_{s\in \mathcal{D}^{(i)}(q)} M_{\alpha^+}(s^{(i)}, s), \qquad  s^{(i+1)} = s^{(i)} + s_q^{(i)}. 
    \label{iterates constraint}
\end{equation}
Then, $m_3(s^{(i+1)})$ satisfies
\begin{equation}
m_3(s^{(i+1)}) - m_{*, \mathcal{D}^{(i)}(q)}\le \left(1-\mu \right)( m_3(s^{(i)}) - m_{*, \mathcal{D}^{(i)}(q)}),
\label{convex linear result}
\end{equation}
where $m_{*, \mathcal{D}^{(i)}(q)} := \min_{s\in \mathcal{D}^{(i)}(q)} m_3(s)$ and $\mu$ satisfying \eqref{cond on mu} is an iteration-independent constant. 

\end{theorem}

\begin{proof} 
With parameters $\iota > 0$, $\kappa > 0$ and $\tilde{p} = 0$, $M_\mu$ is convex with respect to $\mu > 0$. Also, using the condition $H_i \succ 0$, we derive that $M_\mu(s^{(i)}, s)$ is jointly convex for $s \in \R^n$ and $\mu>0$. Therefore, 
\begin{equation} 
M_{\mu}^*(\mu) = \min_{s \in \mathcal{D}^{(i)}(q)} M_\mu(s^{(i)}, s) \qquad \text{is convex for $\mu>0$.}
\label{M alpha tau 1}
\end{equation}
Note that with parameters $\iota > 0$, $\kappa > 0$ and $\tilde{p} = 0$
\begin{equation} 
M_{\mu}^*(1)  = \min_{s \in \mathcal{D}^{(i)}(q)}M_{\alpha^-}(s^{(i)}, s) \le \min_{s \in \mathcal{D}^{(i)}(q)} M(s^{(i)}, s) = m_{*, \mathcal{D}^{(i)}(q)}.
\label{m alpha 1}
\end{equation}
On the other hand, by continuity, we have
\begin{equation} 
M_{\mu}^*(0) = \min_{s \in \mathcal{D}^{(i)}(q)} \lim_{\mu \rightarrow 0 }M_\mu(s^{(i)}, s) =  m_3(s^{(i)}).
\label{m alpha 0}
\end{equation}
Hence,
\begin{eqnarray*}
  m_3(s^{(i+1)})  &=& m_3(s^{(i)}+s_q^{(i)}) = M(s^{(i)}, s_q^{(i)})  
  \\ &\underset{\text{Thm } \ref{thm: Ada Nesterov convex}\text{(2)}}{\le}&     M_{\alpha^+}(s^{(i)}, s_q^{(i)}) = \min_{s\in \mathcal{D}^{(i)}(q)} M_{\alpha^+}({s^{(i)}}, s) \underset{\text{Lemma } \ref{cond alpha convex}}{\le} M_{\mu}({s^{(i)}}, s) 
  \\&\underset{\eqref{M alpha tau 1}}{\le}& 
  \mu M_{\mu}^*(1)  + (1-\mu) M_{\mu}^*(0)
  \underset{\eqref{m alpha 0}, \eqref{m alpha 1}}{\le}
  \mu \min_{s\in  \mathcal{D}^{(i)}(q)} M_{-}({s^{(i)}}, s)  + (1-\mu)   m_3(s^{(i)}) 
 \\&\underset{\text{Thm } \ref{thm: Ada Nesterov convex}\text{(2)}}{\le}&
 \mu \min_{s\in \mathcal{D}^{(i)}(q)}  M({s^{(i)}}, s)  + (1-\mu)   m_3(s^{(i)}) 
= \mu m_{*, \mathcal{D}^{(i)}(q)} + (1-\mu)   m_3(s^{(i)}).
\end{eqnarray*}
Subtracting $m_{*, \mathcal{D}^{(i)}(q)}$ from both sides and rearranging gives the result. 
\end{proof}

{\begin{remark} (Minimizer $s_q^{(i)}$ is well defined)
The quartic regularization term of $M_{\alpha^+}(s^{(i)}, s) \rightarrow \infty$ as $\|s\| \rightarrow \infty$. Thus, there exists a constant $r_q$ such that if $\|s\| > r_q$,  then $M_{\alpha^+}(s^{(i)}, s) \ge f_i = M_{\alpha^+}(s^{(i)}, 0)$. Given that $0 \in \mathcal{D}^{(i)}(q)$, this implies that the minimizer of $M_{\alpha^+}(s^{(i)}, s)$ must satisfy\footnote{A similar argument is used in the proof of Lemma \ref{lemma norm s bound} (see Section \ref{sec: Uniform upper bound} for more details).} $\| s_q^{(i)}\| \le r_q$. Therefore, we can write the constrained optimization problem as
$$ s_q^{(i)} = \argmin_{s\in \mathcal{D}^{(i)}(q)} M_{\alpha^+}(s^{(i)}, s) =  \argmin_{s\in \mathcal{D}^{(i)}(q) \cup \{s \in \R^n,  \|s\| \le r_q\}} M_{\alpha^+}(s^{(i)}, s).$$   
As $\mathcal{D}^{(i)}(q) \cup \{s \in \R^n,  \|s\| \le r_q\}$ is a closed and bounded domain, the minumum of $m_3$ within this domain must be attained, and so $s_q^{(i)}$ is well defined.
\end{remark}}

The proof of Theorem \ref{thm: Ada Nesterov convex} follows a similar approach as Nesterov's proof 
\cite[Theorem 2]{Nesterov2022quartic}. Theorem \ref{thm: Ada Nesterov convex} can be considered as a generalization of the linear convergence result for  Nesterov's method to the case of nonconvex $m_3$ models. In Corollary \ref{corollary: fully convex result}, we present a concise proof illustrating how the linear convergence result of Nesterov's method can be derived from a special case of Theorem \ref{thm: Ada Nesterov convex}.

\begin{corollary} (Linear Rate of Convergence  for Convex $m_3$.) 
Assume $m_3(s)$ is strictly convex for all $s \in \R^n$, and construct $M_{\alpha^+}(s^{(i)}, s)$ as  in Theorem \ref{thm: Ada Nesterov convex} and $s^{(i)}$ is updated by Step 2 of QQR algorithm using \eqref{iterates}. Then, $m_3(s^{(i)})$ satisfies
\begin{equation*}
m_3(s^{(i)}) - m_* \le \left(1-\mu \right)^{i}( m_3(s^{(0)}) - m_*) ,
\end{equation*}
where  $m_* := \min_{s\in \R^n} m_3(s)$ and $\mu$ satisfies $0< \mu \le \min \bigg\{\frac{3\tau - 1}{3\tau + 1}, \left[\frac{1 -2\tau}{1 +2\tau} \right]^{1/3} \bigg\}$ for $\tau \in (\frac{1}{3}, \frac{1}{2})$. 
\label{corollary: fully convex result}
\end{corollary}

\begin{proof}
According to \eqref{nesterov mu bound} and Remark \ref{remark relationship of qqr1 and nesterov}, if $m_3(s)$ is strictly convex, then $M_\mu$ takes the form $M^{\mathcal{N}}_\mu$ with $\tau = \tau_\mathcal{N}$, $q=0$, $\tilde{p}=0$, $\iota = 1$, $\kappa =3$. Additionally, the region $\mathcal{D}^{(i)}(0)$ encompasses the entire space $\mathbb{R}^n$. As a result, the constrained minimization $m_{*, \mathcal{D}^{(i)}(q)} = \min_{s\in \mathcal{D}^{(i)}(q)} m_3(s)$ in \eqref{iterates constraint} simplifies to an unconstrained minimization $m_* = \min_{s\in \R^n} m_3(s)$. Since $m_*$ is independent of the iteration, we can apply induction on $i$ to the inequality $m_3(s^{(i+1)}) - m_* \le (1-\mu)( m_3(s^{(i)}) - m_*)$ and obtain the desired result.
\end{proof}

\begin{remark}
The constant $\mu$ is crucial for determining the linear convergence rate. For nonconvex $m_3$, $\mu$ needs to satisfy \eqref{cond on mu}. In Nesterov's method, which is a special case of the nonconvex framework, we have $q=0$, $\iota=1$, and $\kappa = 3$; therefore, $\tau_{\mathcal{N}}$ is the unique positive root of the equation $\frac{3\tau - 1}{3\tau + 1} = \big[\frac{1 -2\tau}{1 +2\tau} \big]^{1/3}$. This gives $\tau_{\mathcal{N}} = \frac{1}{6}\sqrt{(3+\sqrt{33})} \approx 0.4929$ and $\mu_{\mathcal{N}}\approx 0.193$.
\label{theoretical requirement for mu}
\end{remark}

As we have derived a more general model than the model given by Nesterov, we gain increased flexibility in selecting the parameters $q$, $\iota$, and $\kappa$ to attain a faster linear rate within the region $\mathcal{D}^{(i)}(q)$. In Section \ref{sec: algo for QQRv1}, we will explore the practical implementation of the algorithm and provide detailed discussions on the choices of $q$, $\iota$, and $\kappa$ to achieve a faster linear rate within $\mathcal{D}^{(i)}(q)$. 

It is important to highlight that if the following condition holds
\begin{eqnarray}
\sigma > {\max(s_0^{-2}, 1)} \bigg[\max[-\lambda_{\min}(H), 0] + \hat{L}_H \bigg]>0,
\label{convex sigma bound}
\end{eqnarray}
where $s_0>0$ and $ \hat{L}_H := \max_{v_j \in \R^n, |v_j| = 1, j=1,2,3} \bigg|T_i[v_1][v_2][v_3] \bigg|$, then $m_3(s)$ is convex for all $\|s\|\ge s_0$.  Condition  \eqref{convex sigma bound} implies that for an arbitrarily small $s_0>0$, there exists a sufficiently large $\sigma$ such that $m_3$ is globally convex apart from a small neighbourhood, $\|s\|<s_0$. Several efficient algorithms, such as Nesterov's algorithm \cite{Nesterov2022quartic} or \texttt{QQR-v1}, 
could be employed to minimize a (nearly) convex $m_3$, when initialised in the (convex) region and provided the iterates remain within the respective region of convexity of $m_3$. However, a large $\sigma$ value enforces a stringent step control due to the regularization term, resulting in possibly exceedingly small step sizes for minimizing $m_3$. This will typically be detrimental to the practical performance of the AR$3$ algorithm \cite{birgin2017worst, cartis2020sharp, cartis2020concise}\footnote{It remains to be seen if such a large pre-defined regularization parameter as in \eqref{convex sigma bound} would be detrimental to the theoretical complexity of AR$3$.} 
 when applied to minimizing a general objective $f$ by repeatedly minimizing $m_3$ local models,
as it is crucial then to adapt $\sigma$ to the local landscape, allowing longer steps opportunistically and also to make it independent of difficult-to-calculate or unavailable problem constants. Thus, within each AR$3$ sub-problem, we cannot and should not adjust  $\sigma$ to satisfy the bound in \eqref{convex sigma bound}, and as a result, the sub-problem can become nonconvex when dealing with nonconvex objective functions. This justifies the necessity for QQR methods, which apply to both convex and nonconvex $m_3$.


To conclude this subsection, we note that we did not provide convergence or complexity bounds for nonconvex iterations of \texttt{QQR-v1} in this subsection. The discussions for the convergence of \texttt{QQR-v1} are deferred to Section \ref{sec: theory QQR} and Section \ref{sec: algo for QQRv1}. This is because QQR variant 1 can be considered as a special case of QQR variant 2, as shown in Corollary \ref{link to qqr1}. On the other hand, QQR variant 2, introduced in Section \ref{sec: theory QQR}, contains two adaptive model parameters and is considered a more generic framework. By selecting suitable values for $q$, $\tau$, $\tilde{p}$, $\iota$, and $\kappa$, and adaptively decreasing $\mu$ to a sufficiently small value, QQR variant 1 can conform to the requirements of QQR variant 2 and, thus, exhibit similar convergence behaviour.

\subsection{QQR Variant 2: Two Adaptive Model Parameters}
\label{sec: theory QQR}
In this section, we construct a QQR variant that has two adaptive model parameters, $\alpha_1$ and $\alpha_2$. 
We provide the quadratic and quartic bounds and approximations for $m_3(s)$ at locally convex, locally nonconvex, and locally nearly convex iterations, which we refer to as Cases 1, 2, and 3, respectively.

\begin{tcolorbox}[breakable, enhanced]
 At the $i$th iteration, we iteratively minimize the QQR variant 2 model (\texttt{QQR-v2}) given by:
\begin{equation}
M_\alpha (s^{(i)},s) : = f_i + g_i^Ts+ \frac{\alpha_1}{2} s^T\big[H_i+ \tilde{p} I_n \big] s + \frac{\sigma \alpha_2}{4}\|s\|^4 
\label{qqr v2}
\end{equation}
where $\tilde{p} > 0$ is a fixed constant for locally nearly convex iterations and $\tilde{p} = 0$ otherwise. 
\begin{itemize} \itemsep0em 
    \vspace{-0.2cm}\item  The key idea is to adaptively adjust $\alpha_1$ to switch between locally nonconvex iterations and locally strictly convex iterations, and adaptively adjust $\alpha_2$ to control regularization.
    \vspace{-0.2cm}\item In general, we choose $\alpha_1 \in (0, 1)$ for locally nonconvex iterations and $\alpha_1 > 1$ for locally strictly convex iterations\footnote{For theoretically suitable choices of $(\alpha_1, \alpha_2)$, see Theorems \ref{thm: overall bound} and \ref{thm new bound qqr}; for practical choices of $(\alpha_1, \alpha_2)$, see Section \ref{sec: algo for QQRv2}. }. 
\end{itemize}
\end{tcolorbox}

 We require the following technical lemmas to establish the convergence results for \texttt{QQR-v2}. Proofs of these lemmas will be presented in Section \ref{sec tech proof}.

\begin{lemma} (Upper Bound for $\|s^*\|$) Let $s^* : = \argmin_{s \in \R^n} m_3(s)$ and
$r_c>0$ be the unique implicit solution of the nonlinear equation
$
\sigma = 4\bigg(\frac{\|g\|}{r_c^3} + \frac{\lambda_0}{2r_c^2}+ \frac{\Lambda_0}{6r_c}\bigg),
$
where $ \lambda_0 := \max\{-\lambda_{\min}(H), 0\}$ and $\Lambda_0 := \max_{u\in \R^n, \|u\|=1} {T [u]^3}. $ Then, we have $\|s^*\|\le {r_c}$. Moreover, for any $\|s\|>{r_c}$, $m_3(s) > f_0$. 
\label{lemma norm s bound}
\end{lemma}

\begin{proof}
The proof is given in Section \ref{sec: Uniform upper bound}. 
\end{proof}

\begin{corollary}
There exist iteration-independent positive constants,  such that 
\begin{eqnarray*}
    \max_{u\in \R^n, \|u\|=1, \|s\| \le r_c} \nabla^3 m_3(s) [u]^3  \le L_H, \qquad \max_{u\in \R^n, \|u\|=1, \|s\|\le r_c}   u^T \nabla^2 m_3(s)  u  \le L_g
\end{eqnarray*}
where  $
L_H= \Lambda_0 + 6 \sigma r_c$, $ L_g = |\lambda_{\max}(H)| + \frac{1}{2}\hat{\Lambda}_0 r_c+ 3\sigma r_c^2$,   $\hat{\Lambda}_0 := \max_{v\in \R^n, \|v_1\|=\ldots=\|v_3\|=1} {|T[v_1][v_2][v_3]|}$ and $r_c$ is defined in Lemma \ref{lemma norm s bound}. 
\label{corollary uniform bound}
\end{corollary}

\begin{proof}
    Corollary \ref{corollary uniform bound} is a consequence of  Lemma \ref{lemma tensor update} which is discussed and proved in Section \ref{sec: Lipschitz constant}. Specifically, by substituting $\|s\| \leq  r_c$ in \eqref{LH bound} and \eqref{Lg bound} in Lemma \ref{lemma tensor update} gives the result.  
\end{proof}
We refer to $L_g$ and $L_H$ as the local Lipschitz constants for the gradient and Hessian of $m_3$, respectively (which are global in the neighbourhood $\|s\| \leq  r_c$). A discussion of alternative definitions of Lipschitz constants is given in Remark \ref{lip constant}.

\subsubsection{Bound for Locally Strictly Convex and Nearly Convex Iterations}

In Theorem \ref{thm: overall bound}, we introduce a quadratic quartically-regularised upper/lower bound for $m_3$ that applies in both locally strictly convex and locally nearly convex regimes.

\begin{theorem} Given $s^{(i)} \in \R^{n}$ and {$\|s^{(i)}\| \le r_c$}, let $
M_{\alpha^\pm}(s^{(i)}, s) :  = f_i + g_i^T s + \frac{\alpha_1^\pm}{2} s^T (H_i\pm \tilde p I_n )s + \frac{\sigma}{4}\alpha_2^\pm \|s\|^4. 
$
Let $\alpha_1^\pm = 1\pm a$ and $\alpha_2^\pm = 1\pm d$. Choose $\lambda_{c} > 0$ as a constant, $a$ and $\tilde{p}$ as
\begin{eqnarray}
\textbf{\text{Case 1 (locally strictly convex $\lambda_{\min}[H_i] > \lambda_c)$:}}   && a > 0, \qquad \tilde{p}=0,  \quad   d {>}\frac{L_H^2}{18a\lambda_{c} \sigma},
\label{parameter in nearly convex}
\\
\textbf{\text{Case 2 (locally nearly convex  $\big | \lambda_{\min}[H_i]  \big | \le  \lambda_c)$:}}   && a = 0,\qquad  \tilde{p}>0,  \quad   d {>}\frac{L_H^2}{18  \tilde{p} \sigma},
\label{parameter in convex}
\end{eqnarray}
where $L_H$ is defined in Corollary \ref{corollary uniform bound}. Then, 
\begin{enumerate}
    \item $M(s^{(i)},s)$ satisfies the following upper and lower bounds
$
M_{\alpha^-}(s^{(i)},s) \le M(s^{(i)},s) \le M_{\alpha^+}(s^{(i)},s)
$
for all {$s \in \R^n$} with the equality only at $s \equiv 0$. 
\item  {Letting $
s_+^{(i)} = \argmin_{s\in \R^n} M_{\alpha^+}(s^{(i)}, s) $ and $ s^{(i+1)} := s^{(i)} + s_+^{(i)}$, we have $\|s^{(i+1)}\| \le r_c$.}
\end{enumerate}
\label{thm: overall bound}
\end{theorem}

\begin{proof} 
\textbf{Case 1:} For any $\|s^{(i)}\|\le r_c$,  we have $\big \|T_i  \big\| = \max_{u \in \R^n,  \|u\|=1}  \nabla^3 m_3(s^{(i)})  [u]^3 \le L_H$ by Corollary \ref{corollary uniform bound}. Therefore, 
\begin{eqnarray}
M_{\alpha^+}(s^{(i)},s) - M(s^{(i)},s)  = \frac{a}{2} s^TH_is -  \frac{1}{6}T_i [s]^3 + \frac{d \sigma}{4}  \|s\|^4 \ge  \|s\|^2 \bigg( \frac{a\lambda_{c}}{2}  - \frac{L_H}{6} \|s\|+ \frac{d \sigma}{4} \|s\|^2 \bigg).
 \label{quad2}
\end{eqnarray}
The discriminant of the quadratic equation in \eqref{quad2} is $\Delta_{s} = \frac{L_H^2}{36} - \frac{a\lambda_{c}}{2}  d \sigma.$
Using \eqref{parameter in nearly convex}, we obtain that $\Delta_{s} < 0$ and $\frac{a\lambda_{c}}{2}   - \frac{L_H}{6} \|s\|+ \frac{d \sigma}{4} \|s\|^2>0$. Consequently, $M_{\alpha^+}(s^{(i)},s) - M(s^{(i)},s) \ge 0$ for all {$s \in \R^n$  with equality only at $s \equiv 0 $.} Similar proof follows for $M_{\alpha^-}(s^{(i)},s)$.

Now, we prove the second result. We assume by contradiction that $\|s^{(i)}\| \le r_c$ but  $\|s^{(i+1)}\| =\|s^{(i)}+s_+^{(i)}\| > r_c$. By Result 1), we have
\begin{eqnarray}
M_{\alpha^+}(s^{(i)},s_+^{(i)}) - f_i \ge M(s^{(i)},s_+^{(i)}) - f_i.
\label{bound with eta =1}
\end{eqnarray}
Using $\|s^{(i)}+s_+^{(i)}\| > r_c$, Lemma \ref{lemma norm s bound} and the monotonically decreasing nature of $\{f_i\}_{i >0}$, we have
$$
M(s^{(i)},s_+^{(i)}) - f_i = m_3(s^{(i)}+s_+^{(i)}) - f_i > m_3(s^{(i)}+s_+^{(i)}) - f_0 >0.
$$
Therefore, we obtain that $M_{\alpha^+}(s^{(i)},s_+^{(i)}) - f_i >0$ which is a contradiction to $s_+^{(i)} = \argmin_{s \in \R^n} M_{\alpha^+}(s^{(i)},s) $.

\textbf{Case 2:} We deduce that  $M_{\alpha^+}(s^{(i)},s) - M(s^{(i)},s) \ge \|s\|^2 \bigg( \frac{ \tilde{p}}{2}  - \frac{L_H}{6} \|s\|+ \frac{d \sigma}{4} \|s\|^2 \bigg). $ The same argument follows if we replace $a\lambda_{c}$ by $\tilde{p}$.
\end{proof}

\begin{remark} (A Discussion of Performance Ratio Test) Using the performance ratio to govern progress is a standard numerical technique employed, for example, in AR$p$ algorithms \cite{cartis2007adaptive, cartis2020sharp, cartis2020concise},   trust region algorithms \cite[Sec 3.2.1]{cartis2022evaluation} and the ARC algorithm \cite{cartis2011evaluation}. 
Result 1 of Theorem \ref{thm: overall bound} corresponds to \eqref{ratio test} with $\eta = 1$. Using $M(s^{(i)}, s) \le M_{\alpha^+}(s^{(i)}, s)$, we obtain that 
\begin{eqnarray}
\rho(s) = \frac{m_3(s^{(i)}) - m_3(s^{(i)} + s)}{m_3(s^{(i)}) - M_{\alpha^+}(s^{(i)}, s)} \ge 1 \ge \eta
\label{strong result for ratio}
\end{eqnarray}
for all $\eta \in (0, 1]$ and all $s$. Note that,  for $\eta \in (0, 1]$, Result 2 of Theorem \ref{thm: overall bound} also holds for $s_+^{(i)}$ satisfying $\rho(s_+^{(i)}) \ge \eta$. To see this, we can rewrite \eqref{bound with eta =1} as $M_{\alpha^+}(s^{(i)}, s_+^{(i)}) - f_i \ge \eta \big(M(s^{(i)}, s_+^{(i)}) - f_i\big)$ with $\eta \in (0, 1)$, and then Result 2 of Theorem \ref{thm: overall bound} follows.
\label{eta remark}
\end{remark}

Theorem \ref{thm: overall bound} establishes the upper and lower for $M(s^{(i)}, s)$ for locally convex iterations. With appropriate parameter selections, these bounds can be simplified to \eqref{nesterov upper lower} in the Nesterov model. Corollary \ref{Link to Nesterov Bound} establishes this relationship.

\begin{corollary}
\label{Link to Nesterov Bound}
(Link to Nesterov Bound)  For a strictly convex $m_3$, Theorem \ref{thm: overall bound} simplifies to the upper and lower bounds results in \eqref{nesterov upper lower}, which is given by Nesterov in \cite[Corollary 1]{Nesterov2022quartic}. 
\end{corollary}
\begin{proof}
 To see this, assume that we have $H_i  \succeq \lambda_c I_n \succ  0$ for all $i$, then
\begin{eqnarray}
   0 \prec \nabla^2 m_3(s^{(i)}+s) = H_i  + T_i s+  \sigma (\|s\|^2 + 2 ss^T). 
   \label{convex hessian}
\end{eqnarray}
Applying $s$ twice to \eqref{convex hessian}, we deduce $0 <\lambda_{c} - L_H \|s\| +  3 \sigma \|s\|^2 $. Ensuring the discriminant of the quadratic equation stays negative gives $L_H^2 - 12 \sigma \lambda_c \le 0$. By setting  $a=  \frac{1}{3 \tau}$, we have
$$
d \ge  \frac{L_H^2}{18a\lambda_{c} \sigma} = \frac{L_H^2 \tau}{6\lambda_{c} \sigma}  \ge  2\tau. 
$$ 
where the last inequality uses $L_H^2 \le 12 \sigma \lambda_c$. Therefore, for strictly convex $m_3$, the upper and lower bounds also reduce to $\alpha^{\pm} = {(\alpha_1^{\pm}, \alpha_2^{\pm}) = (1\pm \frac{1}{3 \tau}, 1\pm 2\tau)}. $ Moreover, the restriction for $\alpha_1^{\pm} >0$ and $\alpha_2^{\pm}>0$ gives the bound for $\tau \in (\frac{1}{3}, \frac{1}{2})$.
\end{proof}

\subsubsection{Bound for Locally Nonconvex Iterations}
{
{In a locally nonconvex iteration, where $\lambda_{\min}[H_i] \le -\lambda_c$, instead of proving a upper and lower bound for all $s \in \R^n$, we prove that,  given a constant $\eta \in (0, \frac{1}{3})$,   
there exist suitable values of $(\alpha_1^+, \alpha_2^+)>0$   }
\begin{eqnarray}
\label{nonconvex m alpha}
M_{\alpha^+}(s^{(i)}, s) :  = f_i + g_i^T s + \frac{\alpha_1^+}{2} s^T H_is + \frac{\sigma}{4}\alpha_2^+ \|s\|^4 
\end{eqnarray}
such that  \eqref{ratio test} is satisfied at the minimizer $s_+^{(i)} := \argmin_{s\in \R} M_{\alpha^+}(s^{(i)}, s)$. We first present a technical lemma and then prove the result in Theorem \ref{thm new bound qqr}.

\begin{lemma}
\label{tech lemma 1}
Let  $M_{\alpha^+}$ be defined as in \eqref{nonconvex m alpha} and $s_+^{(i)} := \argmin M_{\alpha^+}(s^{(i)}, s)$. Assume that $\lambda_{\min}[H_i] \le -\lambda_c$ for some $\lambda_c>0$. Then, 
\begin{eqnarray}
g_i^T s_+^{(i)} + \frac{\alpha_1^+}{2} {s_+^{(i)}}^T H_i s_+^{(i)} \le \frac{\alpha_1^+}{2} \lambda_{\min}[H_i] \|{s_+^{(i)}}\|^2 \le  - \frac{\alpha_1^+}{2}\lambda_c \|s_+^{(i)}\|^2 <0. 
\label{eq tech lemma 1}
\end{eqnarray}
\end{lemma}
\begin{proof}
Let $\pm v_1$ be the unit eigenvector that corresponds to the smallest eigenvalue of $H_i$, such that $(\pm v_1)^T H_i(\pm v_1) = \lambda_{\min}[H_i] \le -\lambda_c$. We choose the sign $v_1$ such that $g^T v_1 \le 0$. Since $s_+^{(i)} = \argmin M_{\alpha}(s^{(i)}, s)$, therefore we have
$
    M_{\alpha}(s^{(i)}, s_+^{(i)}) -f_i \le  M_{\alpha}(s^{(i)}, v_1 \|s_+^{(i)}\|) -f_i.
$
Consequently, 
\begin{eqnarray*}
 g_i^T s_+^{(i)} + \frac{\alpha_1^+}{2} {s_+^{(i)}}^T H_i {s_+^{(i)}} +  \frac{\sigma}{4}\alpha_2^+ \|{s_+^{(i)}}\|^4   &\le&   \underbrace{g^T v_1}_{<0}  \|{s_+^{(i)}}\|  + \frac{\alpha_1^+}{2} \underbrace{v_1^T H_iv_1}_{  = \lambda_{\min}[H_i] \le -\lambda_c}  \|{s_+^{(i)}}\|^2 + \frac{\sigma}{4}\alpha_2^+ \underbrace{\|v_1\|^4}_{=1}\|{s_+^{(i)}}\|^4 
\end{eqnarray*}
Thus, we obtain \eqref{eq tech lemma 1}.
\end{proof}

\begin{theorem} Let $\eta \in (0, \frac{1}{3})$ and $M_{\alpha^+}$ be defined as in \eqref{nonconvex m alpha}. Assume that the minimum of $M(s^{(i)},s)$ is not attained at $s =0 $. Given $s^{(i)} \in \R^{n}$ and {$\|s^{(i)}\|< r_c$}, and $\lambda_{\min}[H_i] \le -\lambda_c$. We denote $\alpha_2^+ := \eta(1+ d)$ and choose 
\begin{eqnarray}
\alpha_1^+ \in \bigg[\max \bigg\{\frac{2}{3(1-\eta)}, 1 - \frac{|\lambda_{\min}[H_i]|}{2|\lambda_{\max}[H_i]}\bigg\}, 1\bigg]  \subset (0, 1], \qquad d > \frac{L_H^2}{3 \sigma \lambda_c}.
\label{parameter in nonconvex}
\end{eqnarray}
Then, the following holds
\begin{enumerate}
\item  $M(s^{(i)},{s_+^{(i)}}) -f_i \le \eta \big[M_{\alpha^+}(s^{(i)},{s_+^{(i)}}) -f_i \big]. 
$ Namely, the ratio test in \eqref{ratio test} is satisfied. 
\item  Let $
s_+^{(i)} = \argmin_{s\in \R^n} M_{\alpha^+}(s^{(i)}, s) $ and $ s^{(i+1)} := s^{(i)} + s_+^{(i)}$, then $\|s^{(i+1)}\| \le r_c$.
\end{enumerate}
\label{thm new bound qqr}
\end{theorem}

\begin{proof} 
\textbf{Result 1:} To prove the first result, we define $\mathcal{E}$ as below. We aim to prove that $\mathcal{E} \ge 0$ provided the choice of $(\alpha_1^+, d)$ in Theorem \ref{thm new bound qqr} holds.
\begin{eqnarray*}
  \mathcal{E} :&=&  \eta \bigg[M_{\alpha^+}(s^{(i)},{s_+^{(i)}}) -f_i \bigg] -\bigg[ M(s^{(i)},{s_+^{(i)}}) -f_i     \bigg]=\bigg[f_i-M(s^{(i)},{s_+^{(i)}})\bigg]-\eta\bigg[f_i-M_{\alpha^+}(s^{(i)},{s_+^{(i)}})\bigg]
\\  &=& \eta\bigg[ g_i^T{s_+^{(i)}}  + \frac{\alpha_1^+}{2}{s_+^{(i)}}^T H_i{s_+^{(i)}}   +\frac{\sigma}{4} \big[ \eta^{-1}(1+d)\big]\|{s_+^{(i)}}\|^4\bigg]  -\bigg[ g_i^T{s_+^{(i)}} + \frac{1}{2}{s_+^{(i)}}^T H_i{s_+^{(i)}} + \frac{1}{6}  T_i[{s_+^{(i)}}]^3 + \frac{\sigma}{4} \|{s_+^{(i)}}\|^4\bigg] 
\\   &=& {\underbrace{(\eta-1)}_{<0}g_i^Ts_+^{(i)}} + \frac{1}{2}\underbrace{(\alpha_1^+ \eta-1)}_{<0} {s_+^{(i)}}^T H_i{s_+^{(i)}} - \frac{1}{6}  T_i[{s_+^{(i)}}]^3   + \frac{\sigma}{4} d\|{s_+^{(i)}}\|^4. 
\end{eqnarray*}
Using Lemma \ref{tech lemma 1}, we know that 
\begin{eqnarray}
\underbrace{(\eta-1)}_{<0}g_i^Ts _c + (\eta-1) \frac{\alpha_1^+}{2}{s_+^{(i)}}^T H_i{s_+^{(i)}} \ge \frac{\alpha_1^+}{2} (\eta-1)   \lambda_{\min}[H_i] \|{s_+^{(i)}}\|^2 >0.
\label{temp}
\end{eqnarray}
Substituting \eqref{temp} into $ \mathcal{E}$, we deduce that
\begin{eqnarray}
\notag
  \mathcal{E} & \ge&   \frac{\alpha_1^+}{2} (\eta-1)   \lambda_{\min}[H_i] \|{s_+^{(i)}}\|^2 - (\eta-1) \frac{\alpha_1^+}{2}{s_+^{(i)}}^T H_i{s_+^{(i)}} + \frac{1}{2}\underbrace{(\alpha_1^+ \eta-1)}_{<0} {s_+^{(i)}}^T H_i{s_+^{(i)}}  -\frac{T_i[{s_+^{(i)}}]^3 }{6}   + \frac{\sigma}{4} d\|{s_+^{(i)}}\|^4
  \\ \label{E bound}
  &\ge & \underbrace{\frac{\alpha_1^+}{2}  (\eta-1)   \lambda_{\min}[H_i]}_{> 0}\|{s_+^{(i)}}\|^2 + \frac{1}{2}{\underbrace{(\alpha_1^+ -1)}_{<0} {s_+^{(i)}}^T H_i{s_+^{(i)}} }-\frac{T_i[{s_+^{(i)}}]^3 }{6}   + \frac{\sigma}{4} d\|{s_+^{(i)}}\|^4. 
\end{eqnarray}

\noindent
{\textbf{Case A}}: If ${s_+^{(i)}}^T H_i{s_+^{(i)}} \le 0$,  for all $\alpha_1^+ \in (0, 1]$, following \eqref{E bound}, we have
\begin{eqnarray}
\label{error 2}
 \mathcal{E} 
 \ge \bigg[\overbrace{\underbrace{\frac{\alpha_1^+}{2}  (1-\eta)   \lambda_c }_{> 0} \underbrace{-  \frac{ L_H }{6} \|{s_+^{(i)}}\|}_{<0} +\frac{\sigma}{4} d\|{s_+^{(i)}}\|^2}^{q(\|{s_+^{(i)}}\|)} \bigg] \|{s_+^{(i)}}\|^2
\end{eqnarray}
where we use $T_i[{s_+^{(i)}}]^3 \le L_H \|{s_+^{(i)}}\|$ from Corollary \ref{corollary uniform bound} and
$ \lambda_{\min}[H_i] \le - \lambda_c$. By applying the condition $\alpha_1^+ \ge \frac{2}{3}(1-\eta)^{-1}$ and $d > \frac{L_H^2}{3 \sigma \lambda_c}$, this makes the discriminant of the quadratic equation $q(\|{s_+^{(i)}}\|)$,  $\Delta = \frac{L_H^2}{36} - \frac{\alpha_1^+}{2}  (1-\eta)   \lambda_c \sigma d < 0$. Therefore, we deduce that $ \mathcal{E} \ge 0$. 

\noindent
{\textbf{Case B}}: If ${s_+^{(i)}}^T H_i{s_+^{(i)}} > 0$, let $\alpha_1^+ \in \bigg(\max \bigg\{\frac{2}{3}(1-\eta)^{-1}, 1 - \frac{|\lambda_{\min}[H_i]|}{|2\lambda_{\max}[H_i]|}\bigg\}, 1\bigg] \subset (0, 1]$. 
\begin{itemize}
    \item Since we have $ 1- \frac{-\lambda_{\min}[H_i]}{2\lambda_{\max}[H_i]}  \le \alpha_1^+ \le 1$, we deduce that
$$
0<{\underbrace{(1 -\alpha_1^+ )}_{>0} {s_+^{(i)}}^T H_i{s_+^{(i)}} } \le   \frac{-\lambda_{\min}[H_i]}{2\lambda_{\max}[H_i]} {s_+^{(i)}}^T H_i{s_+^{(i)}} \le \frac{-\lambda_{\min}[H_i]}{2}\|{s_+^{(i)}}\|^2.
$$
where the last inequality uses $\lambda_{\max}[H_i] \ge \frac{{s_+^{(i)}}^T H_i{s_+^{(i)}}}{\|{s_+^{(i)}}\|^2} > 0$ and $\lambda_{\min}[H_i] \le -\lambda_c$.
\end{itemize}
Following \eqref{E bound} and  $T_i[{s_+^{(i)}}]^3 \le L_H \|{s_+^{(i)}}\|$, we have
\begin{eqnarray*}
 \mathcal{E} &\ge&   \bigg[\frac{\alpha_1^+}{2}  (\eta-1)    +\frac{1}{4}\bigg]  \lambda_{\min}[H_i]\|{s_+^{(i)}}\|^2 - \frac{L_H}{6}  \|{s_+^{(i)}}\|^3+ \frac{\sigma}{4} d\|{s_+^{(i)}}\|^4. 
\\
 &\ge&   \bigg[\frac{\alpha_1^+}{2}  (1-\eta)    -\frac{1}{4}\bigg]  \left(-\lambda_{\min}[H_i]\right)\|{s_+^{(i)}}\|^2 - \frac{L_H}{6}  \|{s_+^{(i)}}\|^3+ \frac{\sigma}{4} d\|{s_+^{(i)}}\|^4. 
\end{eqnarray*}
Using $\lambda_{\min}[H_i]\leq -\lambda_c$,
$
 \mathcal{E} \ge   \frac{1}{4}\big[2\alpha_1^+  (1-\eta)    -1\big]  \lambda_c\|{s_+^{(i)}}\|^2 - \frac{L_H}{6}  \|{s_+^{(i)}}\|^3+ \frac{\sigma}{4} d\|{s_+^{(i)}}\|^4. 
$
Using $\alpha_1^+ \ge \frac{2}{3}(1-\eta)^{-1}$,  we arrive at $\frac{\alpha_1^+}{2}(1-\eta) -\frac{1}{4}\ge \frac{1}{3}-\frac{1}{4} = \frac{1}{12}$, 
\begin{eqnarray}
\label{error 3}
 \mathcal{E} \ge   \Bigg[\underbrace{\frac{\lambda_c}{12}}_{>0}  \underbrace{ -\frac{L_H }{6} }_{<0} \|{s_+^{(i)}}\| +\underbrace{\frac{\sigma}{4} d}_{>0}\|{s_+^{(i)}}\|^2\Bigg] \|{s_+^{(i)}}\|^2:= q_2(\|{s_+^{(i)}}\|).
\end{eqnarray}
Using $d > \frac{L_H^2}{3 \sigma\lambda_c}$, the discriminant of $q_2(\|{s_+^{(i)}}\|)$, $\Delta = \frac{L_H^2}{36} - \frac{\lambda_c}{12} \sigma d < 0$. Therefore, we deduce that $ \mathcal{E} \ge 0$. 

In both cases, we have $ \mathcal{E} \ge 0$. Since  the minimum of $M(s^{(i)},s)$ is not attained at $s =0 $, we have $0<M(s^{(i)},{s_+^{(i)}}) -f_i \le \eta \big[M_{\alpha^+}(s^{(i)},{s_+^{(i)}}) -f_i \big] $ and, therefore, \eqref{ratio test} is satisfied. 

\textbf{Result 2:} 
To prove the second result, we assume by contradiction that $\|s^{(i)}\| \le r_c$, but $\|s^{(i+1)}\| =\|s^{(i)}+s_+^{(i)}\| > r_c$. Using Result 1),  Lemma \ref{lemma norm s bound} and the monotonically decreasing nature of $\{f_i\}_{i >0}$, we have
\begin{eqnarray}
M_{\alpha^+}(s^{(i)},s_+^{(i)}) - f_i \ge \eta \big[M(s^{(i)},s_+^{(i)}) - f_i \big] > \eta \big[m_3(s^{(i)}+s_+^{(i)}) - f_0 \big] >0. 
\end{eqnarray}
Therefore, we arrive at $M_{\alpha^+}(s^{(i)},s_+^{(i)}) - f_i >0$ which is a contradiction to $s_+^{(i)} = \argmin_{s \in \R^n} M_{\alpha^+}(s^{(i)},s) $.
\end{proof}

\begin{remark}
    Note that $\alpha_1^+ = 1$ in Theorem \ref{thm new bound qqr} is allowed. This is because we choose $\eta \in (0, \frac{1}{3})$ and therefore then $\frac{2}{3} \le \frac{2}{3}(1-\eta)^{-1} \le 1$. The bound in \eqref{parameter in nonconvex} is well defined. 
\end{remark}

\begin{remark}
Choose $\alpha_1^+ = \max \bigg\{\frac{2}{3}(1-\eta)^{-1}, 1 - \frac{|\lambda_{\min}[H_i]|}{|2\lambda_{\max}[H_i]|}\bigg\}$, in all locally nonconvex iterations, $\alpha_1^+$ has iteration independent upper bound, $\alpha_1^+ \le 1 - \frac{\lambda_c}{2L_g}$.
\label{remark alpha iteration free}
\end{remark}
}

\begin{remark}
Regarding the parameters in the QQR model, we note the following: 
\begin{itemize} \itemsep0em 
\vspace{-0.2cm}\item $\lambda_{c}>0$ is a user-selected fixed positive constant. Usually, we set $\lambda_{c}>\epsilon^{1/3}$ where $\epsilon$ is the first order optimality tolerance. The purpose of $\lambda_{c}$ is to distinguish between the three cases. In Case 2, we can set $\tilde{p}>0$ independently of $\lambda_c$.
\vspace{-0.2cm}\item From  \eqref{parameter in nearly convex},  \eqref{parameter in convex}, \eqref{parameter in nonconvex}, we observe that in locally strictly convex iterations, we have $\alpha_1^+>1$; in locally nearly convex iterations, we have $\alpha_1^+=1$, whereas, in locally nonconvex iterations, we have $0 < \alpha_1^+ \le 1$. This observation explains the choice of $\alpha_1^+$ in Algorithm \ref{TQQR}.
\label{alpha 1 TQQR}
\end{itemize}

\end{remark}

\begin{remark}
$M_{\alpha^{\pm}}$ in Theorems \ref{thm: overall bound} and \ref{thm new bound qqr} preserves local convexity. Note that
$\nabla^2 M_{\alpha^{\pm}}(s^{(i)},s) = \alpha_1^{\pm} H_i + \frac{\sigma}{4} \alpha_2^{\pm} (\|s\|^2 I_n +2ss^T) $ 
where $\alpha_1^{\pm} >0$. If $m_3$ exhibits local convexity with $H_i \succeq 0$, then the bound $M_{\alpha^{\pm}}(s^{(i)},s)$ also shows the same local convexity behaviour as $s \rightarrow \textbf{0}$. The same applies for the locally nonconvex case, $\lambda_{\min}[H_i] \le 0$.  This observation explains Example~\ref{example 1}. 
\label{convexity}
\end{remark}

\begin{remark}
{If we start the \texttt{QQR-v2} Algorithm with $\|s^{(0)}\| \le r_c$ and choose $(a, d, \tilde{p})$ using \eqref{parameter in convex}, \eqref{parameter in nearly convex}, and \eqref{parameter in nonconvex} in each iteration, then Result 1 of Theorems \ref{thm: overall bound} and \ref{thm new bound qqr} indicates that \eqref{ratio test} is always satisfied. Additionally, Result 2 of Theorems \ref{thm: overall bound} and \ref{thm new bound qqr} shows that the iterates are bounded, meaning that the next iterates satisfy $\|s^{(i)}\| \le r_c$.}
\end{remark}

With appropriate parameter selections, the bound in 
Theorems \ref{thm: overall bound} and  \ref{thm new bound qqr} to \eqref{qqr v1} in \texttt{QQR-v1}. Corollary \ref{link to qqr1} provides a direct connection between the bounds derived in  Theorems \ref{thm: overall bound} and \ref{thm new bound qqr} and the specific models of Nesterov and \texttt{QQR-v1}.

\begin{corollary}
(Link to \texttt{QQR-v1} Model) Let the \texttt{QQR-v1} model be defined as in \eqref{qqr v1} with $\kappa, \iota$ and $\tilde{p}$ as specified in Remark \ref{qqr v1 parameter}. Assume that $\|s^{(i)}\|< r_c$ and denote $K_1 := \frac{1}{2} - \frac{1}{6\tau}$ and $K_2 := 1-2(1-q)\tau$, which are both positive constants. If $\mu > 0$ is sufficiently small, satisfying $\mu \le \mu_*$, where
\begin{eqnarray}
   \mu_* = 
\min \bigg\{4^{-1/|\iota|} K_1^{1/\iota} , \bigg(\frac{K_2}{2} \bigg)^{1/\kappa}  , \bigg(\frac{9 \max\{\lambda_c, \tilde{p}\} \sigma }{2L_H^2} {K_2} \bigg)^{1/\kappa} \bigg\}.
\label{mustar}
\end{eqnarray}
Then, $M_\mu(s^{(i)}, s)$ satisfies \eqref{parameter in nearly convex}--\eqref{parameter in convex} in Theorem \ref{thm: overall bound} for locally strictly convex iterations and locally nearly convex iterations.  {$M_\mu(s^{(i)}, s)$ satisfies \eqref{parameter in nonconvex} in Theorem \ref{thm new bound qqr} for locally nonconvex iterations.  
Therefore, the results in Theorem \ref{thm: overall bound} and Theorem \ref{thm new bound qqr} apply for \texttt{QQR-v1}.}  
\label{link to qqr1}
\end{corollary}

\begin{proof}
Let $\alpha_1^{\mu}$ and $\alpha_2^{\mu}$ be defined as in \texttt{QQR-v1} model in \eqref{qqr v1}. We denote \begin{equation}
a = \alpha_1^{\mu} -1 = 2\mu^{-\iota} {K_1} -1, \qquad
d = \alpha_2^{\mu} -1 = \mu^{-\kappa} {K_2} - 1
\end{equation}
where $K_1 = \frac{1}{2} - \frac{1}{6\tau}$ and $K_2 = 1 - 2(1-q)\tau$. 
In \textbf{locally strictly convex iterations}, since $\iota>0$, and $\mu \le (\frac{K_1}{4})^{1/\iota} \le {(\frac{4}{3}K_1)}^{1/\iota} $  gives $a = 2\mu^{-\iota} {K_1} -1 \ge \frac{1}{2} > 0$. We have $a > 0$ which satisfies the requirement for $a$ in local strictly convex iterations in Theorem \ref{thm: overall bound}.
The requirement for $d$ in Theorem \ref{thm: overall bound} is also satisfied since $\mu \le \big(\frac{K_2}{2}\big)^{1/\kappa} $ and $\kappa>0$, then $d = \mu^{-\kappa} {K_2} - 1 \ge \frac{1}{2} \mu^{-\kappa} {K_2} >0$. Thus, we deduce that
\begin{eqnarray}
&& ad =  \bigg| 2\mu^{-\iota} {K_1} -1 \bigg|     \bigg[ \mu^{-\kappa} {K_2} - 1 \bigg] \ge \frac{1}{4} \mu^{-\kappa}{K_2} >\frac{L_H^2}{18 \lambda_c \sigma}
\label{inter step}
\end{eqnarray}
where the last inequality uses the third term in the lower bound of \eqref{mustar}. 
\textbf{Locally nearly convex iterations} follows similarly if we replace $ad$ by $d$ and $\frac{L_H^2}{18 \lambda_c \sigma}$ by $\frac{L_H^2}{18 \tilde{p} \sigma}$ in \eqref{inter step}.

{\textbf{In locally nonconvex iterations,} we have $\alpha^{\mu}_1 = 1$ satisfying \eqref{parameter in nonconvex} in Theorem \ref{thm new bound qqr}. Also, $\alpha^{\mu}_2 = \mu^{-\kappa} \big(1-2(1-q)\tau\big) \ge \frac{2L_H^2}{3 \lambda_c \sigma } \ge \eta( 1 + \frac{L_H^2}{3 \lambda_c \sigma })$  satisfies \eqref{parameter in nonconvex} in Theorem \ref{thm new bound qqr}.}

Therefore, in all cases conditions \eqref{parameter in nearly convex} and \eqref{parameter in convex} are satisfied for Theorem \ref{thm: overall bound}.  
\end{proof}

\subsection{Complexity and Global Convergence of QQR}
\label{Proof of Complexity} 

This section presents a detailed analysis and convergence proofs for the QQR method. The upper bound $M_{\alpha^+}$ in \texttt{QQR-v2} (or $M_{\mu}$ in \texttt{QQR-v1}) plays a crucial role in ensuring the convergence of the QQR algorithm. For simplicity in notation, we will use $M_{\alpha^+}$ to represent the general upper bound in the proofs throughout the subsequent sections. The results apply similarly to $M_{\mu}$. For proving the theoretical results of \texttt{QQR-v2}, we employ the upper bound $M_{\alpha^+}$ throughout this section. Nonetheless, in practical implementations, we minimize $M_{\alpha}$, which is a linear combination of the upper bound $M_{\alpha^+}$ and the lower bound $M_{\alpha^-}$ (see Section \ref{sec: algo for QQRv2}). 

We prove that by iteratively minimizing the upper bound, QQR is guaranteed to converge the local minimum of the quartically-regularised cubic polynomials. Moreover, in cases where the local minimum of $m_3$ is not reached, the function value $m_3(s^{(i)})$ is strictly decreasing with $i$ at each iteration. By imposing slight constraints on the parameters $(a, \tilde{p})$, we demonstrate that the QQR method achieves better complexities than the cubic regularization method in some cases. The complexity bound for the QQR method is determined by $\epsilon$, which is the first-order optimality tolerance, and the parameter $\lambda_c$, which is used to distinguish between different cases. Below is a summary of the complexity bounds for the QQR algorithm. Note that $ \check{c}_1 $ and $ \hat{c}_1 $ are iteration independent constants. 

\begin{itemize}
    \vspace{-0.2cm}\item \textbf{Locally {strictly} convex case ($\lambda_{\min}[H_i] \ge \lambda_c$)} , the error reduction is $m_3(s^{(i)})-  m_3(s^{(i+1)})  \ge  \hat{c}_3 \lambda_c \epsilon,$  or we have a linear convergence in gradient norm.  
    \vspace{-0.2cm}\item \textbf{Locally nonconvex case ($\lambda_{\min}[H_i] \le -\lambda_c$)} , the error reduction is $  m_3(s^{(i)})- m_3(s^{(i+1)})  \ge  \max \big\{ \check{c}_1  \epsilon  \lambda_c^{3/2},   \hat{c}_1  \lambda_c^{3} \big\}$ or we have a constant error reduction. 
    \vspace{-0.2cm}\item \textbf{Locally nearly convex case ($|\lambda_{\min}[H_i] |\le \lambda_c$)} , the error reduction is $  m_3(s^{(i)})- m_3(s^{(i+1)})  \ge  \hat{c}_4 \epsilon^{3/2}$  or we have a linear convergence in gradient norm.  
\end{itemize}
{

We illustrate how the choice of $\lambda_c = \epsilon^\omega$ with $\omega \in [0, \frac{1}{2}]$ affects the complexity bounds in each case. Typically, selecting a larger value of $\lambda_c = \mathcal{O}(1)$ would increase the region of the locally nearly convex case, which has a complexity bound of $\mathcal{O}(\epsilon^{3/2})$. However, this choice provides a better complexity bound for the locally strictly convex and locally nonconvex cases. Conversely, choosing a smaller value of $\lambda_c = \mathcal{O}(\epsilon^{1/2})$ results in a smaller region for the locally nearly convex case (i.e., $|\lambda_{\min}[H_k]|\le \epsilon^{1/2}$) but worsens the complexity bound for the locally strictly convex and locally nonconvex cases.  

\begin{table}[h]
\caption{\small An illustration of the bound in \eqref{function reduction 1}, \eqref{function reduction 2}, \eqref{function reduction 3}  in three convexity cases for different $ \lambda_c = \epsilon^\omega$ with $\omega\in [0, \frac{1}{2}]$.}
\centering
 \begin{center}
 \small{
\begin{tblr}{ |c||c|c|c| } 
 \hline[1.5pt]
\textbf{Order of Bounds} & {$\lambda_{\min}[H_k]\ge \lambda_c = \epsilon^\omega$} & {$ | \lambda_{\min}[H_k]|\le \lambda_c = \epsilon^\omega$} &  {$\lambda_{\min}[H_k] < -\lambda_c = -\epsilon^\omega $}  \\ 
\hline[1pt]
$\omega =0$            &  $\mathcal{O}(\epsilon)$               &  $\mathcal{O}(\epsilon^{\frac{3}{2}})$       & $\mathcal{O}(1)$  
\\ 
$\omega =\frac{1}{4}$  & $\mathcal{O}(\epsilon^{\frac{5}{4}})$     & $\mathcal{O}(\epsilon^{\frac{3}{2}})$  & $\mathcal{O}(\epsilon^{\frac{3}{4}})$
\\ 
$\omega=\frac{1}{3}$   &  $\mathcal{O}(\epsilon^{\frac{4}{3}})$ &  $\mathcal{O}(\epsilon^{\frac{3}{2}})$  & $\mathcal{O}(\epsilon)$   
\\ 
$\omega=\frac{1}{2}$   &  $\mathcal{O}(\epsilon^{\frac{3}{2}})$ &  $\mathcal{O}(\epsilon^{\frac{3}{2}})$  & $\mathcal{O}(\epsilon^{\frac{3}{2}})$   \\ 
 \hline[1.5pt]
\end{tblr}}
\end{center}
\label{table 1}
\end{table}
}

As mentioned in Algorithm \ref{algo1}, we typically set $\lambda_c := \epsilon^{1/3}$, {but note that other values of $\lambda_c = \epsilon^\omega$ with $\omega \in [0, \frac{1}{2}]$ are permissible}. We provide the complexity bounds of the QQR algorithm for $\lambda_c:= \epsilon^{1/3}$. 
\begin{tcolorbox}[breakable, enhanced, title = Summary of Complexity Bounds for QQR]
In iteration $i$, assuming that the minimum of $m_3$ has not been found, i.e., $\|g_i\| > \epsilon$, the QQR algorithm satisfies the following improvements per iteration: 
\begin{enumerate}
    \vspace{-0.2cm}\item \textbf{Locally strictly convex case ($\lambda_{\min}[H_i] \ge \lambda_c$):} $m_3(s^{(i)}) - m_3(s^{(i+1)}) \ge \mathcal{O}(\epsilon^{4/3})$  or we have a linear convergence in gradient norm.  
    \vspace{-0.2cm}\item \textbf{Locally nonconvex case ($\lambda_{\min}[H_i] \le -\lambda_c$):}  $m_3(s^{(i)}) - m_3(s^{(i+1)}) \ge \mathcal{O}(\epsilon)$ or we have a constant error reduction. 
    \vspace{-0.2cm}\item \textbf{Locally nearly convex case ($ |\lambda_{\min}[H_i]| <\lambda_c $):} $m_3(s^{(i)}) - m_3(s^{(i+1)}) \ge \mathcal{O}(\epsilon^{3/2})$  or we have a linear convergence in gradient norm.  
\end{enumerate} 
\end{tcolorbox}

Based on these results, we give the first-order complexity for the QQR method. 

\begin{theorem}
Let $m_3^{\text{low}}$ be the lower bound for $m_3(s)$. Then, there exist positive constants $ \hat{c}_1$,  $\hat{c}_3$, and $\hat{c}_4$ (see Theorems \ref{thm: Convergence and Complexity for Case 3}, \ref{thm: Convergence and Complexity for Case 1}, \ref{thm: Convergence and Complexity for Case 2}) such that, for any $\epsilon \in (0, 1]$, the QQR method (Algorithm \ref{algo1} with $M_{\alpha^+}$ chosen as in Theorems \ref{thm: overall bound} and \ref{thm new bound qqr}) requires at most
    \begin{eqnarray*}
    \#_{\text{QQR}} = \max\big\{\hat{c}_1, \hat{c}_4, \hat{c}_3 \big\}\frac{m_3(s^{(0)}) - m_3^{\text{low}}}{\min\{\epsilon^{4/3}, \epsilon,  \epsilon^{3/2}\}} = \mathcal{O}(\epsilon^{-3/2})
    \end{eqnarray*}
    evaluations of $m_3$ to produce an iterates $s_{\epsilon}$ such that $\|\nabla m_3(s_\epsilon)\| \le \epsilon$. 
\end{theorem}

Although the overall complexity of the QQR method is $\mathcal{O}(\epsilon^{-3/2})$ evaluations, it only represents the worst-case bound among the three cases. We will prove in Corollary \ref{small region} that the $\mathcal{O}(\epsilon^{-3/2})$ case only occurs in a small region where $-\epsilon^{1/3} \le \lambda_{\min}(H_i) \le \epsilon^{1/3}$ and $\epsilon \le \|g_{i+1}\| \le \epsilon^{4/5}$. In our numerical implementation, we have also devised refinement steps to accelerate the convergence to the local minimum when we are close to a local minimum with $\epsilon \leq \|g_{i}\| \leq \epsilon^{4/5}$.
In the other cases, we achieve complexity bounds of order $\mathcal{O}(\epsilon^{-1})$ or $\mathcal{O}(\epsilon^{-4/3})$ evaluations which are better worst-case bounds than that of the cubic regularization method.

\subsubsection{Complexity Bound for Iterations with $\lambda_{\min}[H_i] \le -\lambda_c$}
\label{sec: complexity for nonconvex iteration}
In this subsection, we give proof of the complexity results for the locally nonconvex iterations. 

\begin{theorem}
\label{lemma upper bound m^+} 
Let $s^{(i)}\in \R^n$ be such that $\|g_i\| > \epsilon$ and $\lambda_{\min}[H_i] < - \lambda_c$  with $\lambda_{c}>0$. We construct $M_{\alpha^+}(s^{(i)}, s)$ as defined in Theorem \ref{thm new bound qqr}. 
 $s^{(i)}$ is updated by the QQR algorithm. Assume that  $\|g_{i+1}\| = \| \nabla m_3(s^{(i+1)})\|  > \epsilon$. Then, $m_3(s^{(i)})$ is strictly decreasing with $i$. We have either a constant error reduction $ m_3(s^{(i)})- m_3(s^{(i+1)})  \ge \check{c}_2$ or
\begin{equation}
  m_3(s^{(i)})- m_3(s^{(i+1)})  \ge  {\eta}\max \bigg\{ \check{c}_1  \epsilon  \lambda_c^{3/2},   \hat{c}_1  \lambda_c^{3} \bigg\}. 
   \label{function reduction 1}
\end{equation}
Here, $\hat{c}_1$,  $ \check{c}_1$ and  $\check{c}_2$ are iteration-independent, $\lambda_c$-independent constants and $\epsilon$-independent constants. 
\label{thm: Convergence and Complexity for Case 3}
\end{theorem}

\begin{proof}
 Let $g_v = \frac{g_i}{\|g_i\|}$ be the unit gradient vector. Let $u$ be the unit eigenvector (i.e., $\| u\|=1$) corresponding to $\lambda_{\min}[H_i]$ with the sign of $u$ chosen such that $u^T g_v \le 0$.
Set $c_1>0, c_2 > 0$ as any positive scalars and set $s_1$ as the linear combination of $g_v$ and $u$, such that  $s_1:= c_2 u-c_1g_v$. Substituting $s_1$ into $M_{\alpha^+}$, we have
\begin{eqnarray*}
M_{\alpha^+}(s^{(i)},s_1)  &=&  f_i + {g_i}^T(c_2 u-c_1g_v) +  \frac{\alpha_1^+}{2}  (c_2 u-c_1g_v)^TH_i (c_2 u-c_1g_v)   +  \frac{\sigma}{4}\alpha_2^+\|c_2 u-c_1g_v \|^4  
\\ & \le&  f_i  + c_2\underbrace{g_i^Tu}_{\le 0} -c_1\|g_i\| +
     \frac{\alpha_1^+}{2}\bigg[ c_1^2g_v^TH_ig_v+c_2^2 u^T H_i u-2c_1c_2g_v^TH_iu \bigg]
 +\frac{\sigma}{4}\alpha_2^+(c_1+c_2)^4
\end{eqnarray*}
where the last term comes from $\|c_2 u-c_1g_v \|^4 \le (c_2\|u\|+c_1 \|g_v \|)^4 \le (c_1+c_2)^4$.
Using $c_1>0, c_2>0$, 
\begin{eqnarray*}
u^TH_iu = \lambda_{\min}[H_i] \le -\lambda_{c}, \qquad \text{and}  \qquad
    -g_v^TH_iu = - \underbrace{\lambda_{\min}[H_i]}_{\le 0} \underbrace{g_v^Tu}_{\le 0} \le 0,
\end{eqnarray*} 
we obtain $
M_{\alpha^+}(s^{(i)},s_1)  \le f_i  -c_1\|g_i\| +
     \frac{\alpha_1^+}{2}\bigg[ c_1^2 g_v^T H_i g_v -c_2^2\lambda_{c}\bigg]
 +\frac{\sigma}{4}\alpha_2^+(c_1+c_2)^4.
$
Using $f_i = m_3(s^{(i)})$ and $
g_v^T H_i g_v  \le  L_g$ in  Corollary \ref{corollary uniform bound}, we deduce that
\begin{eqnarray}
m_3(s^{(i)}) - M_{\alpha^+}(s^{(i)},s_1)  
\ge    c_1\|g_i\|-\frac{\alpha_1^+}{2} \big(c_1^2L_g - c_2^2\lambda_{c}\big) - \frac{\sigma}{4}\alpha_2^+ (c_1+c_2)^4.
    \label{nonconvex ineq 1}
\end{eqnarray} 
Also,  for any  $c_1, c_2$, the inequality 
$m_3(s^{(i+1)})  = \min_{s\in \R^n}M_{\alpha^+}(s^{(i)},s) < M_{\alpha^+}(s^{(i)},s_1)$ holds. From Theorem \ref{thm new bound qqr}, we deduce that
\begin{eqnarray}
m_3(s^{(i)}) - m_3(s^{(i+1)})  \ge {\eta} \bigg[ m_3(s^{(i)}) -   M_{\alpha^+}(s^{(i)},s_+^{(i)}) \bigg]\ge {\eta} \bigg[ m_3(s^{(i)}) -   M_{\alpha^+}(s^{(i)},s_1) \bigg]
\label{ratio upper bound}
\end{eqnarray}
where $s_1 := c_2 u-c_1g_v$.
Substituting \eqref{nonconvex ineq 1}  into \eqref{ratio upper bound} gives
\begin{eqnarray*}
{{\eta}}^{-1}\bigg[m_3(s^{(i)}) - m_3(s^{(i+1)}) \bigg] &\ge&  c_1\|g_i\|-\frac{\alpha_1^+}{2} \big(c_1^2L_g - c_2^2\lambda_{c}\big)- \frac{\sigma}{4}\alpha_2^+ (c_1+c_2)^4.
\end{eqnarray*}
Let $
 c_1 := \Bigg[\frac{ \alpha_1^+ L_g } {{\sigma}\alpha_2^+\left[1+\sqrt{2 L_g{\lambda_{c}}^{-1}} \right]^4}\Bigg]^{1/2}>0$ and $  c_2 :=\sqrt{ 2 L_g{\lambda_{c}}^{-1}} c_1>0$, then 
\begin{eqnarray}
{\eta}^{-1}\bigg[m_3(s^{(i)}) - m_3(s^{(i+1)}) \bigg] &\ge&   c_1\|g_i\| + \frac{\alpha_1^+L_g }{2}  c_1^2 -\frac{\sigma \alpha_2^+}{4}\left[1+\sqrt{2 L_g{\lambda_{c}}^{-1}} \right]^4 c_1^4,\notag
\\
&\underset{\|g_i\| > \epsilon}{\ge}&   \epsilon c_1 +\frac{\alpha_1^+L_g}{4} c_1^2 \ge \max \bigg\{ \epsilon c_1,\frac{ \frac{2}{3}(1- \eta)^{-1} L_g }{4}  c_1^2 \bigg\}. 
\label{complexity bound pre}
\end{eqnarray}
The first inequality comes from substituting $c_2$, the second inequality comes from substituting $c_1^2$ in the last term and the last inequality comes from the fact that both  $\epsilon c_1$ and $\frac{\alpha_1^+L_g }{4}  c_1^2$  are positive and  $\alpha_1^+ \ge  \frac{2}{3}(1- \eta)^{-1}$. Therefore, ${\eta}^{-1}\big[m_3(s^{(i)}) - m_3(s^{(i+1)}) \big]$ only needs to be greater than any one of the terms.

\begin{itemize}
    \item 
If $0<\lambda_c < \min \{\frac{L_H^2}{3 \sigma }, 2L_g\} =:\tilde{\lambda}$, 
then, according to \eqref{parameter in nonconvex}, we can set  $d=  \frac{L_H^2}{ \sigma \lambda_c} $  and $\alpha_2^+ = {\eta}^{-1} (1+ d) <2{\eta}^{-1}d \le  \frac{2L_H^2}{ \sigma \lambda_c {\eta}}$. Also, we have $1+ \sqrt{2 L_g{\lambda_{c}}}^{-1} < 2\sqrt{2 L_g{\lambda_{c}}}^{-1}$. Consequently,
$$  
c_1  \ge  \Bigg[\frac{  \frac{2}{3}(1- \eta)^{-1} L_g } { \sigma (\frac{2L_H^2}{ \sigma \lambda_c {\eta}}) \Big(2\sqrt{2 L_g{\lambda_{c}}^{-1}}\Big)^4}\Bigg]^{1/2} = \check{c}_1 \lambda_c^{3/2} 
$$
where $ \check{c}_1$ depends in $\eta$, $L_g$, $L_H$ and $\sigma$, but is iteration independent and $\epsilon$ independent. 
Substituting into \eqref{complexity bound pre}, we arrive at 
$$
  m_3(s^{(i)})- m_3(s^{(i+1)})  \ge   {\eta}  \max \bigg\{ \epsilon  \check{c}_1 \lambda_c^{3/2},\frac{\frac{2}{3}(1- \eta)^{-1}}{4}  \check{c}_1^2 \lambda_c^3\bigg\}.
$$

\item 
If $\lambda_c > \min \{\frac{L_H^2}{3 \sigma }, 2L_g\} =:\tilde{\lambda} >0$,
$
c_1  \ge \bigg[\frac{\frac{2}{3} (1- {\eta} )^{-1}L_g } { \sigma (\frac{2L_H^2}{ \sigma \lambda_c {\eta}})  \big(1+\sqrt{2L_g{\tilde{\lambda}}^{-1}} \big)^4}\bigg]^{1/2} = : \tilde{c}_1 
$ 
which is a constant independent of iteration and $\lambda_c$. We arrive at the constant error reduction
$
  m_3(s^{(i)})- m_3(s^{(i+1)})  \ge  {\eta}  \max \bigg\{ \epsilon  \tilde{c}_1,\frac{ \frac{2}{3}(1- \eta)^{-1} L_g }{4}  \tilde{c}_1^2 \bigg\}:=\check{c}_2.
$
\end{itemize}
\end{proof}

\begin{remark}
\label{remark nonconvex}
Using Theorem \ref{thm: Convergence and Complexity for Case 3}, if assume that $\lambda_c = \epsilon^\omega$ with $\omega \ge 0$, then either there is constant error reduction or
\begin{eqnarray}
  m_3(s^{(i)})-  m_3(s^{(i+1)})  \ge \max \bigg\{  \check{c}_1  \epsilon^{1+\frac{3}{2}\omega} ,  \hat{c}_1\epsilon^{3 \omega} \bigg\} \ge  \mathcal{O}(\epsilon^{\min\{3\omega, 1+\frac{3}{2}\omega\}}).  
  \label{omega complexity bound 1}
\end{eqnarray}
For our choice of $\lambda_c \ge \epsilon^{\frac{1}{3}}$, \eqref{omega complexity bound 1} gives $ m_3(s^{(i)})-  m_3(s^{(i+1)})  \ge  \hat{c}_1 \epsilon.$
\end{remark}


\subsubsection{Complexity Bound for Iterations with $\lambda_{\min}[H_i] \ge \lambda_c$}
In this subsection, we give proof of the complexity results for the locally {strictly} convex iterations. 
To attain the complexity result stated at the beginning of Section \ref{Proof of Complexity}, an additional condition is required when selecting the parameter $a$. Specifically, we need $\epsilon^{1/2} \le a \lambda_c \le \|g_{i}\|^{1/2}$ (or  $\epsilon^{1/2} \le \tilde{p} \le \|g^{(i)}\|^{1/2}$). 

\begin{theorem}
\label{thm: Convergence and Complexity for Case 1}
In iterations with $ \|g_i\|  > \epsilon$ and $ \lambda_{\min}[H_i] \ge \lambda_c$ with $\lambda_{c}>0$, choose $a$ such that
\begin{eqnarray} 
\epsilon^{1/2} \le  a \lambda_c \le  \|g_{i}\|^{1/2}
\label{bound for a}
\end{eqnarray}
and construct $M_{\alpha^+}(s^{(i)}, s)$ as defined in Theorem \ref{thm: overall bound}. The step 
 $s^{(i)}$ is updated by the QQR algorithm. Assume that  $\|g_{i+1}\| = \| \nabla m_3(s^{(i+1)})\|  > \epsilon$. Then, 
\begin{enumerate}
    \vspace{-0.2cm}\item Either there is a linear rate of decrease in gradient norm, such that  $\|g_{i+1}\| \le \rho_{l}\|g_{i}\|$ where $\rho_{l}$ is a fixed constant in $(0, 1)$.
    \vspace{-0.2cm}\item  Or the step size has a lower bound, $ \|s_+^{(i)}\|\ge c_3 \epsilon^{1/2}$. Also, $m_3(s^{(i)})$ is strictly decreasing with $i$ and satisfying
\begin{equation}
  m_3(s^{(i)})-  m_3(s^{(i+1)})  \ge  \hat{c}_3 {\eta} \lambda_c \epsilon.
   \label{function reduction 2}
\end{equation}
Here, $c_3$ and $\hat{c}_3$ are iteration-independent, $\lambda_{c}$-independent and $\epsilon$-independent constants. 
\end{enumerate}
\end{theorem}

\begin{proof}
\textbf{1) Linear Convergence in $\|g_i\|$}:  If $\|g_{i+1}\| \le \rho_{l} \|g_{i}\| $, we have a linear rate of decrease in $\|g_i\|$. 

\noindent
\textbf{2) Function Value Decrease}:  In the hard case, we have
$
    \|g_{i+1}\| >  \rho_{l}\|g_{i}\|. 
$
We write
$$
g_{i+1} = \nabla M(s^{(i)}, {s_+^{(i)}}) = \nabla M(s^{(i)}, {s_+^{(i)}})-\nabla M_{\alpha^+}  (s^{(i)}, {s_+^{(i)}}) + \nabla M_{\alpha^+}(s^{(i)}, {s_+^{(i)}}). $$ By the first-order optimality conditions of $ M_{\alpha^+}$, $\nabla M_{\alpha^+}(s^{(i)}, {s_+^{(i)}}) =0$. Thus, 
\begin{eqnarray}
g_{i+1} =   \nabla M(s^{(i)}, {s_+^{(i)}})-\nabla M_{\alpha^+}  (s^{(i)}, {s_+^{(i)}}) = a H_{i} s_+^{(i)} - \frac{1}{2} T_i[s_+^{(i)}]^2 + d \sigma \|s_+^{(i)}\|^2 s_+^{(i)} . \label{argument 1}
\end{eqnarray}
where $ d = \frac{L_H^2}{18 a \lambda_c}$ (Theorem \ref{thm: overall bound}). Using triangle inequality, we obtain that 
\begin{eqnarray}
\epsilon < \|g_{i+1}\| \le a \lambda_c \|{s_+^{(i)}}\| + \frac{ L_H }{2} \|{s_+^{(i)}}\| ^2 + a^{-1}\lambda_{c}^{-1}  \tilde{d} \|{s_+^{(i)}}\|^3
\label{gi+1 bound}
\end{eqnarray}
where $ \tilde{d} := \frac{L_H^2}{18}$ is an iteration-independent constant and $d = a^{-1} \tilde{d}$. At least one term in the right hand side of \eqref{gi+1 bound} needs to be greater than $\frac{1}{3}\|g_{i+1}\|$. Therefore, one of the following inequalities is true
\begin{eqnarray}
    \|{s_+^{(i)}}\| &\ge & \frac{\|g_{i+1}\| }{3a \lambda_c} \ge \frac{1}{3} \|g_{i+1}\|^{1/2} \frac{\|g_{i+1}\|^{1/2} }{\|g_{i}\|^{1/2}} \underset{ \|g_{i+1}\| >  \rho_{l}\|g_{i}\|} {\ge} \frac{1}{3}\rho_{l}^{1/2}  \|g_{i+1}\|^{1/2},
    \label{s bound 1}
    \\ \text{or}, \qquad  
    \|{s_+^{(i)}}\| &\ge & \bigg[\frac{2\|g_{i+1}\| }{3L_H}\bigg]^{1/2} \ge \big(\frac{2 }{3L_H}\big)^{1/2} \|g_{i+1}\|^{1/2},
    \\ \text{or}, \qquad  
    \|{s_+^{(i)}}\| &\ge & \bigg[\frac{ a\lambda_{c}\|g_{i+1}\| }{3\tilde{d}}\bigg]^{1/3} \underset{\eqref{bound for a}} {\ge} \bigg[\frac{ \epsilon^{1/2}\|g_{i+1}\| }{3\tilde{d}}\bigg]^{1/3} \ge(3\tilde{d})^{-1/3} {\epsilon^{1/6}}\|g_{i+1}\|^{1/3} .
\label{s bound 3}
\end{eqnarray}
By setting $c_3 :=\min \big\{\frac{1}{3}\rho_{l}^{1/2}, \big(\frac{2 }{3L_H}\big)^{1/2}, (\frac{L_H^2}{6})^{-1/3} \big\} $, using $\|g_{i+1}\| > \epsilon$, we proved the lower bound for step size $ \|s_+^{(i)}\|\ge c_3 \epsilon^{1/2}$. Note that $c_3$ is an iteration independent and $\lambda_c$ independent constant.  Using \eqref{bound for a} gives $a >0$,  $M(s^{(i)},s) \le M_{\alpha^+}(s^{(i)},s)$ in Theorem \ref{thm: overall bound} and {Remark \ref{eta remark}}, we deduce that
\begin{eqnarray}
     {\eta^{-1}} \bigg[ m_3(s^{(i)}) - m_3(s^{(i+1)}) \bigg] &=&  m_3(s^{(i)}) - M(s^{(i)}, s_+^{(i)} ) \ge  m_3(s^{(i)}) 
 - M_{\alpha^+}(s^{(i)}, s_+^{(i)} ) 
 \label{argument 2}
 \\  &=& - g_i^T {s_+^{(i)} }^T  - \frac{\alpha_1^+}{2} {s_+^{(i)} }^T H_i s_+^{(i)}   - \frac{\alpha_2^+ }{4} \|s_+^{(i)} \|^4. 
 \label{argument 5}
\end{eqnarray}
The first-order optimality conditions of $ M_{\alpha^+}$ gives $
   -g_i =  \alpha_1^+ H_is_+^{(i)}   + \alpha_2^+ \sigma \big\|s_+^{(i)} \big\|^2 s_+^{(i)}. $
Consequently, 
\begin{eqnarray}
   {\eta^{-1}} \bigg[ m_3(s^{(i)}) - m_3(s^{(i+1)}) \bigg] \ge  \frac{\alpha_1^+}{2}{s_+^{(i)}}^T H_i  s_+^{(i)}   + \frac{3\alpha_2^+}{4} \sigma  \|s_+^{(i)} \|^4 \ge  \frac{\alpha_1^+}{2}  \lambda_c   \|s_+^{(i)}\|^2 \ge \frac{\lambda_c}{2}  \|s_+^{(i)}\|^2
     \label{argument 4}
\end{eqnarray}
where we use $\alpha_1^+ = 1+a>1$, $\alpha_2^+ \ge 0 $ and $\lambda_{\min} [H_i] \ge \lambda_c>0$. Finally, using $ \|s_+^{(i)}\|\ge c_3 \epsilon^{1/2}$, we obtain 
$$
     m_3(s^{(i)}) - m_3(s^{(i+1)}) \ge \frac{\lambda_c}{2}c_3^2  {\epsilon}  = \lambda_c
     {\hat{c}_3\epsilon} 
$$
where 
$
\hat{c}_3 :=  \frac{1}{2}c_3^2 =  \frac{1}{2}\min \big\{\frac{\rho_l}{9}, \frac{2}{3L_H}, (\frac{L_H^2}{6})^{-2/3} \big\}
$ is  an iteration-independent, $\lambda_c$ iteration-independent and $\epsilon$-independent constant.
\end{proof}

\begin{remark}
Since $ \|g_{i}\| \ge \epsilon$, the bound for $a$ in \eqref{bound for a} is well-defined. 
\end{remark}

\begin{remark}
\label{remark convex}
Using Theorem \ref{thm: Convergence and Complexity for Case 2}, if assume that $\lambda_c = \epsilon^\omega$ with $\omega \ge 0$, then either there is a linear rate of decrease in gradient norm or
\begin{eqnarray}
  m_3(s^{(i)})-m_3(s^{(i+1)})  \ge  \hat{c}_3 \epsilon^{1+\omega}.
  \label{omega complexity bound 2}
\end{eqnarray}
For our choice of $\lambda_c \ge \epsilon^{\frac{1}{3}}$, \eqref{omega complexity bound 2} gives $ m_3(s^{(i)})-  m_3(s^{(i+1)})  \ge \hat{c}_3  \epsilon^{4/3}.$
\end{remark}

\begin{remark}  (Extension to Locally Nearly Convex Case)
In locally nearly convex iterations, if we impose a similar condition on $\tilde{p}$ as stated in \eqref{bound for a}, specifically, $\epsilon^{1/2} \le \tilde{p} \le \|g^{(i)}\|^{1/2}$, then under the iteration where $H_i+\tilde{p}I_n \succ 0$, the analysis presented in \eqref{argument 1}--\eqref{argument 4} follows in a similar manner. 
\end{remark}


\subsubsection{Overall Complexity Bound}
\label{sec: complexity for nearly iteration}

We have calculated complexity bounds for locally {strictly} convex and nonconvex iterations in the cases where $\lambda_c > \epsilon^{1/3}$ (as discussed in Remark \ref{remark nonconvex} and Remark \ref{remark convex}). The only case that is left to address is when $0 \le \lambda_c \le \epsilon^{1/3}$.  To address this, Theorem \ref{thm: Convergence and Complexity for Case 2} will demonstrate that the function value reduction in this scenario is at least $\mathcal{O}(\epsilon^{3/2})$. Importantly, Theorem \ref{thm: Convergence and Complexity for Case 2} is universally applicable for all $\lambda_c > 0$ and all convexity cases. Theorem \ref{thm: Convergence and Complexity for Case 2} establishes that the worst-case function value decrease achieved by the QQR method is at least as good as that observed in locally convex iterations.

\begin{theorem}
\label{thm: Convergence and Complexity for Case 2}
In locally strictly convex or nonconvex iterations with $ \|g_i\|  > \epsilon$ and $ |\lambda_{\min}[H_i]| \ge \lambda_c$, choose $a$ such that
$$
\epsilon^{1/2} \le  |a \lambda_c| \le  \|g_{i}\|^{1/2}. 
$$
In locally nearly convex iterations with $ \|g_i\|  > \epsilon$ and $ |\lambda_{\min}[H_i]| \le \lambda_c$, choose $a$ such that
$$
\epsilon^{1/2} \le  \tilde{p} \le  \|g_{i}\|^{1/2}. 
$$
We construct $M_{\alpha^+}(s^{(i)}, s)$ as defined in Theorems \ref{thm: overall bound} and \ref{thm new bound qqr} and update $s^{(i)}$ by the QQR algorithm. Assume that $\|g_{i+1}\| = \| \nabla m_3(s^{(i+1)})\|  > \epsilon$; then, either  $\|g_{i+1}\| \le \rho_{l}\|g_{i}\|$ where $\rho_{l}$ is a fixed constant in $(0, 1)$. Or, $\|s_+^{(i)}\|\ge c_3 \epsilon^{1/2}$, and $m_3(s^{(i)})$ is strictly decreasing with $i$ and satisfies
\begin{equation}
  m_3(s^{(i)})-  m_3(s^{(i+1)})  \ge \hat{c}_4 \epsilon^{2/3}\|g_{i+1}\|^{5/6}  \ge \hat{c}_4 \epsilon^{3/2}. 
   \label{function reduction 3}
\end{equation}
where $c_3$ and $\hat{c}_4$ are iteration-independent, $\lambda_c$ iteration-independent and $\epsilon$-independent constants.
\end{theorem}
\begin{proof}
We give the proof for the case of $ \lambda_{\min}[H_i] \ge \lambda_c$, and the proof for locally nearly convex cases follows similarly if we replace $a \lambda_c$ by $\tilde{p}$. This proof is applicable for all $\lambda_c>0$.  

If $\|g_{i+1}\| \le \rho_{l} \|g_{i}\| $, we have a linear rate of decrease in $\|g_i\|$. Otherwise, if $\|g_{i+1}\| \ge \rho_{l} \|g_{i}\| $, by a similar deduction as in \eqref{argument 1}--\eqref{s bound 3}, we arrive at the bound for $\|s_+^{(i)}\|$ in \eqref{s bound 1}--\eqref{s bound 3}. According to Theorem 8.2.8~\cite{cartis2022evaluation}, the  first-- and second--order optimality condition of $ M_{\alpha^+}$ gives
\begin{eqnarray*}
-g_i =  \alpha_1^+ H_is_+^{(i)}   + \alpha_2^+ \sigma \|s_+^{(i)} \|^2 s_+^{(i)},\qquad \alpha_1^+ H_i  \succeq - \alpha_2^+  \sigma  \|{s_+^{(i)}} \|^2I_n.
\end{eqnarray*} 
Using these optimality conditions, from \eqref{argument 5}, we deduce that
\begin{eqnarray}
     m_3(s^{(i)}) - m_3(s^{(i+1)}) \ge  \frac{\alpha_1^+}{2}{s_+^{(i)}}^T H_i  s_+^{(i)}   + \frac{3\alpha_2^+}{4} \sigma  \|s_+^{(i)} \|^4 \ge \frac{\alpha_2^+}{4} \sigma  \|s_+^{(i)} \|^4 \ge \frac{L_H^2}{72} \|g_i\|^{-1/2} \|s_+^{(i)} \|^4
     \label{argument general}
\end{eqnarray}
where the last inequality comes from  $\alpha_2^+ = 1+d >d = \frac{L_H^2}{18 a \lambda_c}$, and  $\alpha_2^+ \sigma \ge \frac{L_H^2}{18|a| \lambda_c } \ge \frac{L_H^2}{18} \|g_i\|^{-1/2}$. Substituting \eqref{s bound 1}--\eqref{s bound 3} into \eqref{argument general},
\begin{eqnarray*}
     m_3(s^{(i)}) - m_3(s^{(i+1)}) &\ge&  \frac{L_H^2}{72}  \min \bigg\{\frac{\rho_{l}^2 }{3} \frac{\|g_{i+1}\|^{2}}{ \|g_i\|^{1/2}}, \frac{4}{9L_H^2} \frac{\|g_{i+1}\|^{2}}{ \|g_i\|^{1/2}},  \big(\frac{L_H^2}{6}\big)^{-4/3} {\epsilon^{2/3}} \frac{\|g_{i+1}\|^{4/3}}{\|g_i\|^{1/2}}\bigg\}
     \\ & \underset{ \|g_{i+1}\| >  \rho_{l}\|g_{i}\|} {\ge}&  \min \bigg\{\frac{L_H^2\rho_{l}^{5/2} }{216}\|g_{i+1}\|^{3/2}, \frac{\rho_{l}^{1/2}}{162}\|g_{i+1}\|^{3/2}, \frac{\rho_{l}^{1/2}}{6L_H^{3/2}}\epsilon^{2/3}\|g_{i+1}\|^{5/6}\bigg\}
     \\&\ge& \hat{c}_4 \epsilon^{2/3}\|g_{i+1}\|^{5/6}  \ge \hat{c}_4 \epsilon^{3/2} 
\end{eqnarray*}
where $\hat{c}_4 :=  \min \bigg\{\frac{L_H^2\rho_{l}^{5/2} }{216}, \frac{\rho_{l}^{1/2}}{162}, \frac{\rho_{l}^{1/2}}{6L_H^{3/2}}\bigg\}$ is an iteration independent and $\lambda_c$ independent constant. 
\end{proof}

Our complexity proofs for  Theorem \ref{thm: Convergence and Complexity for Case 1} and Theorem \ref{thm: Convergence and Complexity for Case 2} follow a similar structure as the complexity proof for ARC. These results were first obtained by Nesterov and Polyak \cite[Thm 3.3.5]{Nesterov2006cubic}. Other works can be found in \cite{birgin2017worst, cartis2020concise}. To conduct the complexity analysis, we employ a typical construction that involves demonstrating lower bounds on the model decrease and the length of the step at each iteration. The connection between the gradient of model decrease and the gradient of function value decrease is established using \eqref{argument 4}. {While our proof assumes exact minimization of $M_{\alpha^+}$ such that $\nabla M_{\alpha^+}(s^{(i)}, s^{(i)}_+) = 0$, the same order of complexity bounds as in Theorem \ref{thm: Convergence and Complexity for Case 1} and in Theorem \ref{thm: Convergence and Complexity for Case 2} 
can be derived when considering approximate minimization with a step-termination condition, such as $\|\nabla M_{\alpha^+}(s^{(i)}_+)\| \le \theta \|s^{(i)}_+\|^3$ for some $\theta \in (0, 1)$. This step-termination condition is also utilized in \cite{cartis2020concise}, and alternative variants exist \cite{cartis2011adaptive, gratton2022adaptive}. 
The proof incorporates standard techniques found in \cite[Lemma 2.3]{birgin2017worst} and in \cite[Lemma 3.3]{cartis2020concise}, which demonstrate that the step cannot be arbitrarily small compared to the criticality conditions. We delegate the use of approximate termination conditions in the subproblem solution to future work and to current numerical results.} 
\begin{corollary}
Using Theorem \ref{thm: Convergence and Complexity for Case 3}--\ref{thm: Convergence and Complexity for Case 1}, let $\lambda_c = \epsilon^{1/3}$. Then, only in the small region where $-\epsilon^{1/3} \le \lambda_{\min}(H_i) \le \epsilon^{1/3}$ \textbf{and} $\epsilon\le \|g_{i+1}\| \le \epsilon^{4/5}$, there is a linear rate of decrease in gradient norm, or 
\begin{equation*}
  m_3(s^{(i)})-  m_3(s^{(i+1)}) \ge  \hat{c}_4 \epsilon^{3/2}.
\end{equation*}
Otherwise, we have a linear rate of decrease in gradient norm, constant error reduction or $
  m_3(s^{(i)})-  m_3(s^{(i+1)}) \ge \min\{\hat{c}_1,\hat{c}_3, \hat{c}_4 \}\epsilon^{4/3}. 
$ Constants are defined in  Theorem \ref{thm: Convergence and Complexity for Case 3}--\ref{thm: Convergence and Complexity for Case 1} and are iteration-independent, $\lambda_c$ iteration-independent and $\epsilon$-independent constants.
\label{small region}
\end{corollary}

\begin{proof}
 For $\lambda_c = \epsilon^{1/3}$,  if $|\lambda_{\min}(H_i)|>\lambda_c \ge \epsilon^{1/3}$  according to  Remark \ref{remark nonconvex} and Remark \ref{remark convex}, we have error reduction of at least $\mathcal{O}(\epsilon)$  (or constant decrease) for the locally nonconvex case and $\mathcal{O}(\epsilon^{4/3})$ (or linear decrease in gradient) for locally {strictly} convex case.  If $\|g_{i+1}\| \ge \epsilon^{4/5}$, then according to Theorem \ref{thm: Convergence and Complexity for Case 2}, we have either linear rate linear decrease in gradient or
 $  m_3(s^{(i)})-  m_3(s^{(i+1)}) \ge \hat{c}_4 \epsilon^{\frac{2}{3}}\epsilon^{\frac{4}{5} \frac{5}{6}}  \ge \hat{c}_4 \epsilon^{4/3}.$
\end{proof}

In the implementation of the QQR algorithm, when we are close to a local minimum with $\epsilon \leq \|g_{i}\| \leq \epsilon^{4/5}$, we often employ a refinement step to accelerate the convergence to the local minimum. Specifically, we set a small value for the {quartically regularization parameter} and relax the control over the regularization to enhance asymptotic convergence. As a result, in our numerical experiments, we rarely encounter the worst-case scenario where the convergence rate is $\mathcal{O}(\epsilon^{3/2})$. A more comprehensive discussion on the implementation of QQR and the refinement step can be found in Section \ref{sec_implementation}.

To conclude this section, it is important to highlight that there may be several alternative techniques for constructing the quadratic quartic local models within the context of the QQR method. For instance,  as illustrated in \cite{zhuquartic}, we can use the gradient information, i.e., the (local) gradient Lipschitz constant $L_g$, to obtain the local upper/lower bounds for $m_3(s)$ at $s = s^{(i)}$. Let $\tau_\pm>0$, then the gradient-based bounds are
$$M_{\alpha^\pm}(s^{(i)},s) = m_3({s^{(i)}}) + g_i^Ts+  s^T\bigg[\left(\frac{1}{2}-  \frac{1}{3\tau_\pm}\right)\nabla^2_s m_3(s^{(i)})  \pm \frac{L_g}{3\tau_\pm}I_n\bigg]s  + \left(1\pm 2 \tau_\pm\right) \frac{\sigma}{4} \|s\|^4 
$$
where  $\tau_+>1$ and $0 < \tau_- < \frac{1}{2}$. The gradient-based upper bounds also have a comparable structure to Nesterov's method. By selecting specific values for $\tau_\pm$, the gradient-based upper bounds can be simplified to match the expression in Nesterov's method. 
To avoid complicating the presentation of our general QQR method, we will not delve into the specifics of a gradient-bound-based QQR method. Instead, we focus on creating a QQR method that exploits (local) second and third-order
smoothness, adapting the parameters automatically during the iterations. 

\section{Bounded Minimizer and Lipschitz Constants}
\label{sec tech proof}

{
In this section, we prove the existence of a uniform upper bound, $r_c$, such that for $s^* = \argmin_{s \in \R^n} m_3(s)$, we have $\|s^*\| \le r_c$. We discuss local Lipschitz constants and provide technical lemmas to prove  Corollary \ref{corollary uniform bound}.
}

\subsection{Proof of Lemma \ref{lemma norm s bound}}
\label{sec: Uniform upper bound}

\begin{proof}
We define a function $K(r): \R^+ \rightarrow \R$, 
$
K(r):= \frac{\|g\|}{r^3} + \frac{\lambda_0}{2r^2}+ \frac{\Lambda_0}{6r}.
$
$K(r)$ is a monotonically decreasing function with $r>0$. $K(r) \rightarrow +\infty$ as $r\rightarrow 0$ and $K(r) \rightarrow 0$ as $r\rightarrow +\infty$. Thus, we can find a unique constant $r_c>0$ that satisfies
$\frac{\sigma }{4}= \frac{\|g\|}{r_c^3} + \frac{\lambda_0}{2r_c^2}+ \frac{\Lambda_0}{6r_c}
$. Therefore, we re-write $m_3(s)$ as, 
$
m_3(s) = f_0 + g^T s +  \frac{1}{2} s^T H s+ \frac{1}{6}T[s]^3+ \bigg[\frac{\|g\|}{r_c^3} + \frac{\lambda_0}{2r_c^2}+ \frac{\Lambda_0}{6r_c}\bigg]\|s\|^4. 
$
Let  $\xi_s := \frac{\|s^*\|}{{r_c}}$. Assume for contradiction the global minimum of $m_3(s)$ happens at $\|s^*\|>{r_c}$ (namely $\xi_{s^*}>1$).  Then, 
\begin{eqnarray}
m_3(s^*) = f_0+ \underbrace{\bigg[g^Ts^*+ \xi_s^3 \|g\|\|s^*\|\bigg]}_{>0} + \frac{1}{2} \underbrace{ {s^*}^T\bigg[H+\lambda_0 \xi_s^3 I_n\bigg]s^*} _{>0} + \frac{1}{6}\underbrace{\bigg[T [s^*]^3 +\Lambda_0 \xi_s^3 \|s^*\|^3\bigg] }_{>0}. 
\label{mid step 2}
\end{eqnarray}
Thus, $m_3(s^*) > f_0.$ This contradicts the hypothesis that $s^*$ is the global minimum. 
\end{proof}

\subsection{Local Lipschitz Constants}
\label{sec: Lipschitz constant}

We prove that technical Lemma for Corollary \ref{corollary uniform bound}. 

\begin{lemma}
For any $s \in \R^n$, $ L_H(s)  $ satisfies
\begin{equation}
 L_H(s) : =   \max_{u\in \R^n, \|u\|=1}  \nabla^3 m_3(s) [u]^3  \le  \Lambda_0 + 6 \sigma\|s\|
 \label{LH bound}
\end{equation} 
where $\Lambda_0 := \max_{u\in \R^n, \|u\|=1} {T [u]^3}$. Similarly,  
\begin{equation}
L_g(s) : =   \max_{u\in \R^n, \|u\|=1}   u^T \nabla^2 m_3(s)  u   \le \max \bigg\{\lambda_{\max}(H) + \frac{1}{2} \hat{\Lambda}_0\|s\|+ 3\sigma \|s\|^2, 0 \bigg\}
\label{Lg bound}
\end{equation} 
where $\hat{\Lambda}_0 := \max_{v\in \R^n, \|v_1\|=\|v_2\|=\|v_3\|=1} {|T[v_1][v_2][v_3]|}$.
\label{lemma tensor update}
\end{lemma}
\begin{proof}
For any $s = [s_1, \dotsc, s_n]^T\in \R^n$, the third order derivative is $
 \nabla^3 m_3(s) = T +\frac{\sigma}{4}\nabla^3 (\|s\|^4)
$
where $\nabla^3 (\|s\|^4) \in \R^{n\times n \times n}$ is a symmetric tensor with
\begin{eqnarray}
\text{the $j$th face of $\nabla^3 (\|s\|^4)$} = 8(se_j^T + e_js^T+s_j I)
\end{eqnarray}
for $ 1\le j \le n$. $e_j \in \R^{n}$ is the unit vector with one on the $j$th entry and zero otherwise and $I_n \in \R^{n\times n}$ is the identity matrix. 
We first prove the bound for $L_H(s)$. For all $u \in \R^n$ with $\|u\|=1$, 
\begin{eqnarray*} 
\nabla^3 (\|s\|^4)[u]^3 = 8\sum_{j=1}^n u_j u^T \bigg[se_j^T + e_js^T+s_j I_n\bigg] u = 8 \sum_{j=1}^n (2u_j^2 u^Ts + u_js_j \|u\|^2) = 24 u^Ts.  
\end{eqnarray*}
Thus, $\nabla^3 m_3(s) [u^3]= T[u^3] +\frac{\sigma}{4}\nabla^3 (\|s\|^4)  \le  \Lambda_0  + 6 \sigma u^Ts.$ Using Cauchy-Schwarz inequality $u^Ts \le \|s\|\|u\| \le \|s\|$ gives \eqref{LH bound}. Similarly for $L_g(s)$, we have $$u^ T\nabla^2 m_3(s) u = u^ T H u + \frac{1}{2}T[s][u]^2 + \sigma \bigg[\|s\|^2+ (u^Ts)^2\bigg].$$  Therefore, by Cauchy-Schwarz inequality and triangle inequality $\max_{u\in \R^n, \|u\|=1}  \big| u^T \nabla^2 m_3(s) 
 u\big| \le \lambda_{\max}(H) + \frac{1}{2} \hat{\Lambda}_0 \|s\|+ 3\sigma \|s\|^2. $
\end{proof}

\begin{remark}
\label{lip constant}
(Alternative Definitions of Lipschitz Constant)
Let $m_3$ be defined as in \eqref{m3} and the third order Taylor model of $m_3(s)$ at $s = s^{(i)}$ is
$
    T_3(s^{(i)}, s)=  f_i + g_i^T s +  \frac{1}{2} s^TH_i s  + \frac{1}{6}  T_i  [s]^3. 
$
Note that $T_3(0, s)=  f_0 +  g^T s +  \frac{1}{2} s^TH s  + \frac{1}{6}  T [s]^3$. 
Let $\hat{L}_H(s^{(i)})$ be the Lipschitz constant of $\nabla^2 T_3(s)$, such that
\begin{eqnarray*}
\big \|\nabla^2 T_3(s^{(i)}, s_1)-\nabla^2 T_3(s^{(i)}, s_2) \big \|_{[2]} = \big \|   T_i  (s_1 - s_2)  \big \|_{[2]}
\le \hat{L}_H(s^{(i)}) \|s_1-s_2\|
\end{eqnarray*}
for all $s_1, s_2\in \R^n$. $\|\cdot\|_{[p]}$ is the tensor norm recursively induced by the Euclidean norm on the space of $p$th-order tensors 
$\|\mathcal{T}\|_{[p]}:= \max_{\|v_1\| =  \dotsc= \|v_p\| = 1\|} \big|\mathcal{T}[v_1, \dotsc, v_p] \big|. $ Then, 
\begin{equation*}
     \hat{L}_H(s^{(i)})  = \big \|T_i  \big\|_{[3]} = \max_{v \in \R^n,  \|v_1\| =  \dotsc= \|v_3\| = 1}  \big|T_i[v_1][v_2][v_3] \big|.  
\end{equation*}
$L_H(s^{(i)})$ is the maximum absolute value of $T_i[v_1][v_2][v_3]$ by applying any three unit vectors. But in the complexity and convergence proof, we can use a tighter tensor norm for $T_i$ which is the maximum absolute value of $T_i[u]^3$ by applying any same unit vectors three times. For any $s^{(i)} \in \R^n$, we define
\begin{eqnarray*}
 L_H(s^{(i)})  : = \big \|T_i  \big\|_{[3]^*} = \max_{u \in \R^n,  \|u\|=1}  T_i [u]^3 .  
\end{eqnarray*}
Note that $
 L_H(s^{(i)}) \le   \hat{L}_H(s^{(i)}) $
for all $s^{(i)} \in \R^n$.  Note that we can choose the sign of $\pm u \in \R^n$ such that $T_i [u]^3 >0$. Therefore, by definition, $L_H(s^{(i)}) \ge 0$ for all $s^{(i)} \in \R^n$.
\end{remark}

\section{Implementations of the QQR Algorithmic Framework}
\label{sec_implementation}

In this section, we introduce two implementations of the QQR algorithmic framework, leveraging the theoretical results established in the previous section. The first implementation, which is built on \texttt{QQR-v1} in Section \ref{sec: locally strictly convex Iterations}, involves the adaptive adjustment of a single parameter. The second implementation, built on \texttt{QQR-v2} in Section \ref{sec: theory QQR}, incorporates the adaptive adjustment of two parameters. The numerical set-up and examples are provided in  Section \ref{sec: Numerics example}. It is important to note that we do not explicitly compute the (local) Lipschitz constants $L_g$ or $L_H$ that appear in the theoretical analysis. Instead, we employ adaptive parameter adjustment(s) strategies to achieve the desired theoretical bounds. {Note that \texttt{QQR-v1} and \texttt{QQR-v2} (Algorithms \ref{alg: qqr v1}--\ref{alg: qqr v2}) are constructed within the framework of Algorithm \ref{algo1}, with more specific parameter selections.}

\subsection{Implementation of Variant 1}
\label{sec: algo for QQRv1}

\textbf{Choice of $q$:} We have demonstrated that for sufficiently small $\mu$ satisfying \eqref{cond on mu}, the iterates exhibit linear convergence with a contraction factor of $(1-\mu)$ within the region $\mathcal{D}^{(i)}(q)$ (Theorem \ref{thm: convex complexity result}). 
Also, the choice of $q$ within the range of $0$ to $3$ determines the trade-off between the size of the linear convergence region $\mathcal{D}(q)$ and the impact of the regularization term in $M_{\mu}$ (Theorem \ref{thm: Ada Nesterov convex}). Setting $q = 0$ expands the linear convergence region $\mathcal{D}(0)$, but at the cost of smaller step sizes due to the dominant regularization term. Conversely, as $q$ approaches $3$, the regularization term becomes less influential, allowing for larger step sizes, but there is a risk of the iterates falling outside the linear convergence region $\mathcal{D}(3)$. Therefore, selecting a value of $q$ that strikes a balance between these considerations is crucial. In line with the global optimality condition stated in Theorem \ref{thm: arc}, which requires $\mathcal{D}^{(i)}(q)$ to satisfy $H_i + T_i[s] + \sigma \|s\|^2I_n \succeq 0$, we opt for $q=2$.

\textbf{Choice of $\iota$ and $\kappa$:} In the case of locally convex iterations, we set $\iota = \kappa = 1$, which ensures that $M_{\mu}$ represents a linear combination of the upper bound and the lower bound. 

\textbf{Rate of Linear Convergence $\tau_*$ under Convexity}: The optimal linear rate is determined by the largest possible $\mu$ that satisfies the condition in \eqref{cond on mu}. Substituting  $\iota = \kappa =1 $ and $q = 2$, we have
$    \mu_c = \max_{\tau \in (\frac{1}{3}, 1)}\min \big\{\frac{3\tau - 1}{3\tau + 1},\frac{1 -\frac{2}{3}\tau}{1 +\frac{2}{3}\tau} \big\}. $
This gives a unique optimal root, $\mu_c = \frac{1}{7}(11 - 6 \sqrt{2}) \approx 0.359$ and consequently, the optimal (local) linear rate $\tau_* = \frac{\sqrt{2}}{2}$.

\begin{remark}
Compared to Nesterov's methods,  \texttt{QQR-v1}  achieves a faster (local) linear rate of convergence. The contraction factor for Nesterov's methods is $1 - \mu_{\mathcal{N}} \approx 0.807$, whereas the contraction factor for  \texttt{QQR-v1}  is $1 - \mu_c \approx 0.641$. 
\label{mu selection}
\end{remark}

\textbf{Number of Unsuccessful Iterations:} We typically initialise $\mu$ with $\mu = \mu_0$ and then adaptively decrease $\mu$ by updating $\mu: = \eta_0 \mu$, where $0 < \eta_0 < 1$.  If the iteration is locally strictly convex, we set $\iota = 1$. We have proven that choosing $\mu < \mu_c \approx 0.359$ ensures (local) linear convergence. Therefore, it takes at most $\lceil \log(\mu_c / \mu_0) / \log(\eta_0) \rceil + 1$ unsuccessful iterations to decrease $\mu$ below $\mu_c$ before we can accept the step.
If the iteration is locally nonconvex, according to Corollary \ref{link to qqr1}, we require $0 < \mu < \mu_*.$ Therefore, it takes at most $\lceil \log(\mu_* / \mu_0) / \log(\eta_0) \rceil + 1$ unsuccessful iterations to decrease $\mu$ below $\mu_*$ before we can accept the step.
If the iteration is nearly convex, then we add a perturbation to the second-order term by setting $\tilde{p} >0$.

\begin{algorithm}[!ht]
\caption{\small \texttt{QQR-v1} Algorithm}\label{alg: qqr v1}
\small{\textbf{Initialization}: $s^{(0)} = 0$, $i=0$, $0 < \rho_1 < {\frac{1}{3}}< \rho_2$, $0 < \gamma_0 < 1<\gamma_1$. 
\\Set $\epsilon > 0$, usually $\epsilon  =10^{-5} \sim 10^{-8}$, 
$\lambda_{c} = \epsilon^{1/3}$. 
\\Use the model in \eqref{qqr v1} (with initialization $\kappa=\iota = 1$, $q=2$, $\mu = \frac{1}{2}$, $\tau_* =\frac{\sqrt{2}}{2}$ and $\tilde{p}=0$). 

\While{$\|g_i\|  >\epsilon $}
{
Compute ${s_+^{(i)}} := \argmin_{s\in \R^n} M_\mu(s^{(i)}, s)$ and $\rho(s_+^{(i)})$ as in \eqref{ratio test}. 

\If{$\rho(s_+^{(i)})\ge \rho_1$, \underline{Successful}}
{Accept the iterate. 
$s^{(i+1)} := s^{(i)}+{s_+^{(i)}}, i=i+1, \tilde{p}=0$. 
\\ \If{$\rho(s_+^{(i)})\ge \rho_2$ \underline{{Very Successful}}} 
{\textbf{Increase $\mu$} by $\mu = \gamma_1 \mu$}
}

\If{$\rho(s_+^{(i)})< \rho_1$, \underline{Unsuccessful}} 
{Reject the iterate.
Evaluate $\lambda^{(i)} = {\lambda_{\min}[H_i]}$.}

\eIf{$|\lambda^{(i)}| \le \lambda_c$}
{Add perturbation $\Tilde{p} >0$. }
{
\textbf{Decrease $\mu$} by $\mu = \gamma_0 \mu$. If $\lambda^{(i)} < -\lambda_c$, set {$\alpha_1^{\mu} := 1$.} 
}}}
\end{algorithm}

\subsection{Implementation of Variant 2}
\label{sec: algo for QQRv2}

In the implementation of \texttt{QQR-v2}, we initialise $M_\alpha$ with $\alpha_1 = \alpha_2 = 1$ and  $\tilde{p}=0$.  The complete algorithm for \texttt{QQR-v2} is provided in Algorithm \ref{alg: qqr v2}. {We note that both the smallest eigenvalue $\lambda_{\min}[H_i]$ and largest one $\lambda_{\max}[H_i]$ are needed in the running of the algorithm.}

\begin{algorithm}[!ht]
\caption{\small \texttt{QQR-v2} Algorithm}
\label{alg: qqr v2}
\small{\textbf{Initialization}: $s^{(0)} = 0$, $i=0$, $0 < \rho_1  < {\frac{1}{3}} < \rho_2$, , $0 < \eta_0 < 1<\eta_1$, $\gamma_2>1$.  
\\Set $\epsilon  > 0$, usually $\epsilon  =10^{-5} \sim 10^{-8}$, 
$\lambda_{c} = \epsilon^{1/3}$. 

Use the model in \eqref{qqr v2} (with initialization $\alpha_1=\alpha_2=1$ and $\tilde{p}=0$): 
\\
\While{$\|g_i\|  >\epsilon$}
{
Compute ${s_+^{(i)}} := \argmin_{s\in \R^n} M_\alpha(s^{(i)}, s)$ and $\rho(s_+^{(i)})$ as in \eqref{ratio test}. 

\If{$\rho(s_+^{(i)})\ge \rho_1$, \underline{Successful}} 
{Accept the iterate. $s^{(i+1)} := s^{(i)}+s_+^{(i)}, i=i+1, \tilde{p}=0$. 
\\ \If{$\rho(s_+^{(i)})\ge \rho_2$, \underline{Very Successful}} 
{\textbf{Decrease $\alpha_2$} by $\mu = \eta_0 \alpha_2$}
}
\If{$\rho(s_+^{(i)}) < \rho_1$, \underline{Unsuccessful}} 
{Reject the iterate. Evaluate $\lambda^{(i)} = {\lambda_{\min}[H_i]}$.
\\
\eIf{$|\lambda^{(i)}| \le \lambda_c$}
{Add perturbation $\Tilde{p}> 0$. }
{Increase $\alpha_2$ by $\alpha_2 :=\eta_1 \alpha_2$, and
\\ If $\lambda^{(i)}< -\lambda_{c}$, set {$\alpha_1 := \max \bigg\{\frac{2}{3}(1-\eta)^{-1}, 1 - \frac{|\lambda_{\min}[H_i]|}{2|\lambda_{\max}[H_i]}\bigg\}$}, 
\\ If $\lambda^{(i)} > \lambda_c$, increase $\alpha_1 := \gamma_2 \alpha_1$. 
}}}}
\end{algorithm}

\begin{remark} (Other Variations of the Algorithm)
Note that when $s^{(i)}$ is close to a local minimum of $m_3(s)$, setting the regularization parameter to $\min \big\{\frac{\sigma}{4}, \|g_i\| \big\}$ is a standard numerical technique \cite{cartis2007adaptive} to speed up the asymptotic convergence in the convex region. As the iterates approach the local minimum, the coefficient of the quartic-order term tends to $\epsilon$, allowing for more relaxed control over the regularization. Other numerical techniques that speed up asymptotic convergence, such as the Newton method, can also be applied.
\end{remark}

\subsection{Numerical Set-up and Illustrations}
\label{sec: Numerics example}
We construct quartically-regularised cubic polynomials to test the QQR algorithm's performance. Specifically, we define $m_{3}(s) =f_0 + g^Ts + \frac{1}{2}  s^T H s + \frac{1}{6} T [s]^3+ \frac{\sigma}{4}  \|s\|^4$, where $f_0=0$, and the coefficients $g$, $H$, and $T$ are generated as follows: 
\begin{eqnarray*}
\quad g = \texttt{a*randn(n, 1)}, \quad  H = \texttt{b*symm(randn(n, n))},  \quad  T = \texttt{c*symm(randn(n, n, n))}.
\end{eqnarray*}
$g$ is a random $n$-dimensional vector multiplied by a scaling factor \texttt{a}; $H$ is a random symmetric $n\times n$ matrix multiplied by \texttt{b}; and $T$ is a random symmetric $n\times n\times n$ tensor multiplied by \texttt{c}. Here, \texttt{n} represents the dimension of the problem ($s\in \R^n$), and \texttt{symm(randn())} denotes a symmetric matrix or tensor whose entries follow a normal distribution with mean zero and variance one. The parameters \texttt{a}, \texttt{b}, \texttt{c}, and $\sigma$ are chosen differently to test the algorithm's performance under various scenarios. We set the stopping criteria to be $\|g_i\|  < \epsilon$ where $\epsilon= 10^{-5}$ unless otherwise specified. The \textit{iteration counts} record the number of times we minimize the quadratic model with {quartic} regularization. The {iteration counts} include both successful and unsuccessful iterations. The \textit{functions'/derivatives' evaluation counts} record the number of times we evaluate the function and derivative values. The {evaluation counts} only include the successful iterations.

\begin{figure}[!ht]
    \centering
    \includegraphics[width = 10cm]{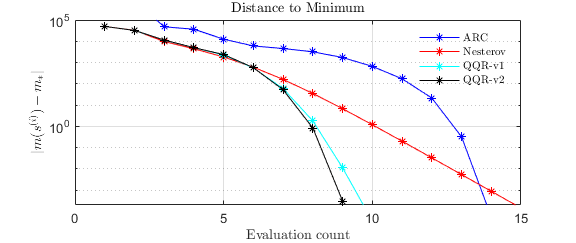} 
    \caption{\small The parameters are \texttt{a} $=80$, \texttt{b} $=80$, \texttt{c} $=80$, and $\sigma =80$, $n=100$.}
    \label{Nesterov vs qqr1} 
\end{figure}

According to Figures \ref{Nesterov vs qqr1},  \texttt{QQR-v1}  performs competitively with Nesterov's method in terms of both evaluation counts. In locally strictly convex iterations (which start from the third evaluation count), there are two main factors contributing to the rate of  \texttt{QQR-v1}. Firstly, as illustrated in Remark \ref{mu selection},  \texttt{QQR-v1}  has a contraction factor of $0.641$, while Nesterov's method has a contraction factor of $0.807$. In addition to  \texttt{QQR-v1}  having a faster theoretical linear rate, the model we minimize in each step of \texttt{QQR-v1} also incorporates a linear combination of (local) upper bound $M_{\alpha^+}$ and (local) lower bound $M_{\alpha^-}$. In contrast, in Nesterov's method, the model we minimize in each step is only the (local) upper bound. In locally nonconvex iterations, neither method guarantees a linear rate of convergence. 
Among these methods,  \texttt{QQR-v2}  requires the least number of function/derivative evaluation counts to find the local minimum of $m_3(s)$. The ARC method only adjusts the parameter in the regularization term, whereas \texttt{QQR-v2} adjusts both the coefficients in the second-order term and the regularization term.

\section{Numerical Results: Comparison Across Methods}
\label{sec_numerics}
This section provides a comprehensive comparison of minimization algorithms for the AR$3$ sub-problem discussed in this paper. These algorithms include the ARC method, Nesterov's method (for convex $m_3$ only), \texttt{QQR-v1} method, and \texttt{QQR-v2} method. The ARC method was proposed by Cartis et al. in \cite{cartis2011adaptive}, while Nesterov's method was introduced in \cite{Nesterov2022quartic}. \texttt{QQR-v1} and \texttt{QQR-v2} were proposed in this paper.

To assess the robustness of the minimization solvers, we devised multiple sets of AR$3$ sub-problems across different dimensions, $n$, ranging from $5$ to $100$. These test sets comprise AR$3$ sub-problems that feature varying gradients, Hessians, tensors, and regularization terms. Specifically, we consider singular, ill-conditioned, concave or convex, diagonal, and indefinite Hessians, as well as large and small tensor terms, dense tensor terms, and ill-conditioned tensor terms. 

Table \ref{table convex H} to Table \ref{table large tensor} report the numerical results for these test sets. The coefficients in the test sets are generated as described in Section \ref{sec: Numerics example}. The results presented in each table are obtained over $10$ random problems randomized using MATLAB functions \texttt{rand()} or \texttt{randn()}. The \textit{relative minimum values} of a method is obtained by dividing the minimum value achieved by this algorithm by the mean of all the minimum values obtained by all algorithms. In all of the test cases, none of the algorithms resulted in unsuccessful cases, as they all converged to a local minimum that satisfies $\|g_i\| < \epsilon$ and $\lambda_{\min}[H_i] > -\epsilon^{1/2}$ with $\epsilon= 10^{-5}$.  The best-performing results for each problem set are highlighted.

The numerical results show that overall, our proposed methods, including  \texttt{QQR-v1}  and   \texttt{QQR-v2},  consistently perform competitively with the ARC method and  Nesterov's method in terms of the number of evaluations of functions/derivatives or iterations required. Our proposed methods are robust enough to handle ill-conditioning in the Hessian term (Table \ref{table ill con H}) and perform well for the minimization of $m_3$ with singular, concave, diagonal, or indefinite Hessians  (Table \ref{table convex H}, \ref{table ill con H}). Additionally, QQR is effective in finding local minima for quartic polynomials with ill-conditioned tensor terms, diagonal tensor terms, and large and directional tensors Table \ref{table large tensor}. Last but not least, our proposed algorithm can locate the local minimum for a wide range of positive $\sigma$ values (Table \ref{change sigma}). This enables us to dynamically adjust $\sigma$ while utilizing QQR for the AR$3$ method in minimizing the objective function $f(x)$. We provide some highlights specific to each method and each problem sets. 

\begin{itemize} \itemsep0em 
\vspace{-0.2cm}\item Our proposed \texttt{QQR-v1} method is a generalization of Nesterov's method and performs competitively with it in terms of function/derivative evaluations and iterations for all dimensions. A notable finding is that \texttt{QQR-v1} performs competitively with Nesterov's method even for convex quartically-regularised cubic polynomials (Table \ref{table convex H}). The reason behind this is that our implementation of \texttt{QQR-v1} uses a model that is a linear combination of both the local upper bound and lower bound, while Nesterov's method only utilizes the upper bound. When dealing with convex and locally strictly convex $m_3$, our method can converge to the local minimum with just a few iterations and function evaluations.

\vspace{-0.2cm}\item Our proposed \texttt{QQR-v2} method generally requires the least function/derivative to locate the local minimum of $m_3$ among all methods. All methods have comparable CPU times, with \texttt{QQR-v2} requiring slightly less time than the other methods.

\vspace{-0.2cm}\item QQR can minimize ill-conditioned problems (Table \ref{table ill con H}). Interestingly, in the ill-conditioned test case in Table \ref{table ill con H},  QQR locates the local minimum of $m_3$ with very few function/derivative evaluations. This can be attributed to the Hessian being diagonal and positive semi-definite, which results in faster convergence, as seen in the convex quartically-regularised cubic polynomial minimization (Table \ref{table convex H}). 


\vspace{-0.2cm}\item Based on our numerical experiments, the nearly convex case happens very rarely in the QQR iterates, even for ill-conditioned Hessian terms or singular Hessian terms. This size of the perturbation term $\tilde{p}$ does not significantly affect the convergence behaviour in practice. 
\end{itemize}

The selection of which method to use for minimizing the quartically regularised cubic polynomial is contingent on the intended objective of the optimization. If the goal is to satisfy the complexity requirements outlined in \cite{cartis2020concise} by attaining a local minimum for the AR$3$ sub-problem, then  \texttt{QQR-v2}  is generally favored. 

\begin{table}[!htbp]
\centering
\caption{\small \textbf{Convex $m_3$, Convex $H$ or Concave $H$:}  \texttt{a} $=80$, and $\sigma =80$.  In the \textbf{first table}, $m_3$ is convex, $H = \texttt{symm(30*(randn(n)+neye(n)))}$ and \texttt{c} $=1$; In the \textbf{second table}, $m_3$ is nonconvex but locally convex at $s=0$, $H = \texttt{symm(30*(randn(n)+neye(n)))}$ and \texttt{c} $=80$. In the \textbf{third table}, $m_3$ is concave, $H = \texttt{symm(30(randn(n)-n*eye(n)))}$  and \texttt{c} $=80$. The CPU time is calculated as the average CPU time of the solver for solving all problems across all dimensions.  $\rho_1 = 0.3, \rho_2 = 3, \eta_0 = 0.5, \eta_1 = 2, \gamma_0 = 0.8, \gamma_1 = 1.2, \gamma_2 = 1.1$ for Tables \ref{table convex H} to \ref{table large tensor}.}
    \label{table convex H}
{\scriptsize
\begin{tblr}{
  colspec = {c|rrr|rrr|rrr|r},
}
\hline[2pt]
   & \SetCell[c=3]{}Relative Min. & & & \SetCell[c=3]{}Iterations. & & &  \SetCell[c=3]{}Evaluations.& & & CPU time\\
\hline[1pt]
Method  &   $n=5$ &  50 &  100 & 5 &  50 &  100  &   5 &  50 &  100 & (s)\\
\hline[0.5pt]
    \textbf{\color{red}ARC}  &  -1 & -1 & -1
    &  5 & 4 & 3
    &   5 & 4 & 3& 0.024\\
    \textbf{Nesterov} &  -1 & -1 & -1
    &  12.9 & 14 & 15 
    &  12.9 & 14 & 15 & 0.077\\
   \textbf{\color{blue}QQR-v1}  &  -1 & -1 & -1
   & 6 & 7.3 & 7
    &6 & 7 & 7 & 0.049\\
   \textbf{\color{blue}QQR-v2} &  -1 & -1 & -1
    &  \bb{2.9} & \bb{2} & \bb{2} 
    &  \bb{4} & \bb{3} & \bb{3} & \bb{0.018} \\
\hline[0.5pt]
\hline[0.5pt]
    \textbf{\color{red}ARC}  &  -1 & -1 & -1
    &  5 & 4 & 4  
    &   5 & 4 & 4 & 0.032 \\
   \textbf{\color{blue}QQR-v1}  &  -1 & -1 & -1
   &6.5 & 7 & 7  &6.5 & 7 & 7 & 0.054\\
   \textbf{\color{blue}QQR-v2} &  -1 & -1 & -1
    &  \bb{4} & \bb{3} & \bb{3} 
    &  \bb{4} & \bb{3} & \bb{3} & \bb{0.029} \\
\hline[0.5pt]
\hline[0.5pt]
    \textbf{\color{red}ARC}  &  -0.998 & -0.995 & -0.997
    &  11.7 & 24.1 & 27.6& 
   8 & 14.4& 21.4 & 0.129\\
   \textbf{\color{blue}QQR-v1}  &   1 & 1 & 1
   &9.7 & 22& 23.4     &9.6 & 20.7&22.6 &  0.145\\
   \textbf{\color{blue}QQR-v2} &  -0.998 & -0.995 & -0.997
    &  \bb{7.4} & \bb{12.2} & \bb{17} 
    &  \bb{7.3} & \bb{11.8} & \bb{17}  & \bb{0.106}\\
\hline[2pt]
\end{tblr}
}
\end{table}

\begin{table}[!htbp]
\centering
\caption{\small \textbf{Ill-Conditioned $H$, Diagonal $H$ or Singular $H$:} \texttt{a} $=80$, \texttt{c} $=80$ and $\sigma = 80$. $H$ is a diagonal ill-conditioned matrix with diagonal entries uniformly distributed in $[0, 10^{10}]$.   }
\label{table ill con H}
{\scriptsize
\begin{tblr}{
  colspec = {c|rrr|rrr|rrr|r},
}
\hline[2pt]
   & \SetCell[c=3]{}Relative Min. & & & \SetCell[c=3]{}Iterations. & & &  \SetCell[c=3]{}Evaluations.& & & CPU time \\
\hline[1pt]
Method  &   $n=5$ &  50 &  100 & 5 &  50 &  100  &   5 &  50 &  100 & (s)\\
\hline[0.5pt]
    \textbf{\color{red}ARC}  &  -1 & -1 & -1 
    &  5.3 & 6.5 & 5.3 
      &  5.3 & 6.2 & 5.3 & 0.042 \\
   \textbf{\color{blue}QQR-v1}   &  -1 & -1 & -1 
   & 7.3 & 7.5& 9.7
   & 7 & 7.4& 9.2 & 0.062\\
   \textbf{\color{blue}QQR-v2}   &  -1 & -1 & -1 
   &\bb{4.1} & \bb{4.8} & \bb{4.7}
   &\bb{4.1} & \bb{4.8} & \bb{4.5} & \bb{0.037} \\
\hline[2pt]
\end{tblr}
}
\end{table}

\begin{table}[!htbp]
\centering
    \caption{\small \textbf{Changing $\sigma$: }  From top to bottom, \texttt{a} $=80$, \texttt{b} $=80$,  \texttt{c} $=80$ and $\sigma =5, 300$, respectively.  The experiments show the algorithm's ability to handle a range of $\sigma$ values. Generally, smaller values of $\sigma$ for $m_3(s)$ require more function/derivative evaluations for convergence across all algorithms.}
     \label{change sigma}
{\scriptsize
\begin{tblr}{
  colspec = {c|rrr|rrr|rrr|r},
}
\hline[2pt]
   & \SetCell[c=3]{}Relative Min. & & & \SetCell[c=3]{}Iterations. & & &  \SetCell[c=3]{}Evaluations.& & & CPU time \\
\hline[1pt]
Method  &   $n=5$ &  50 &  100 & 5 &  50 &  100  &   5 &  50 &  100 & (s)\\
\hline[0.5pt]
    \textbf{\color{red}ARC}  &  -1 & -1 &-1
    &  \bb{11.8} & 24.8 & 28
    & \bb{8.1} & 15.2 &20.5& 0.12\\
   \textbf{\color{blue}QQR-v1}  & -0.84 & -0.9958 & -0.897 
   &31.9&19.2 & 37.9
   &11.3 & 18.5 &26.7 & 0.19\\
   \textbf{\color{blue}QQR-v2} & -1  & -0.9958  & -0.9978
    &  16.1 & \bb{12.9} & \bb{18.1} 
    &  14.5 & \bb{12.8} & \bb{17.6} & \bb{0.10} \\
\hline[0.5pt]
\hline[0.5pt]
    \textbf{\color{red}ARC}  &  -1 & -1 & -1
    & 7 & 9.8 & 20.4
    &6.7 &8.8 & 12.4 & 0.09\\
   \textbf{\color{blue}QQR-v1}   &  -1 & -1 & -1
   &7.6& 11.3 & 14 &
  7.6& 10.6& 13.4 & 0.10\\
   \textbf{\color{blue}QQR-v2} &  -1 & -1 & -1
    &\bb{4.8} & \bb{7}& \bb{14}
    &  \bb{4.8} & \bb{7} & \bb{10.2} & \bb{0.08}\\
\hline[2pt]
\end{tblr}
}

\end{table}

\begin{table}[!htbp]
\centering
    \caption{\small \textbf{Small Tensor Term v.s. Large Tensor Term: }  From top to bottom, \texttt{a} $=80$, \texttt{b} $=80$, $\sigma =80$, and \texttt{c} $=300, 10$ respectively.}
        \label{table small tensor}
{\scriptsize
\begin{tblr}{
  colspec = {c|rrr|rrr|rrr|r},
}
\hline[2pt]
   & \SetCell[c=3]{}Relative Min. & & & \SetCell[c=3]{}Iterations. & & &  \SetCell[c=3]{}Evaluations.& & & CPU time \\
\hline[1pt]
Method  &   $n=5$ &  50 &  100 & 5 &  50 &  100  &   5 &  50 &  100 & (s)\\
\hline[0.5pt]
    \textbf{\color{red}ARC}  &  -1 & -1 & -0.979
    &  \bb{9.5} & 19.3 & 26.9
    &  \bb{8} & 15.7 & 19 & 0.128 \\
   \textbf{\color{blue}QQR-v1}   &  -0.84 & -1 & -1
   & 39.5 & 18.5 & 20.7
   & 16.9 & 17.8 & 20.4 & 0.163\\
   \textbf{\color{blue}QQR-v2}    &  -1 & -1 & -1
    &  12.8 & \bb{12.5} & \bb{16.9}
    &  11.8 & \bb{12.5} & \bb{16.7} & \bb{0.123} \\
\hline[0.5pt]
\hline[0.5pt]
    \textbf{\color{red}ARC}  &  -1 & -1 & -1 
    &  7.1 & 11.4 & 15.1
    &  6.3 & 7.9 & 9.9 & 0.0695 \\
   \textbf{\color{blue}QQR-v1}   &  -1 & -1 & -1 
   & 6.8 & 8.6 & 10
   & 6.7 & 8.6 & 9.8& 0.0691 \\
   \textbf{\color{blue}QQR-v2} &  -1 & -1 & -1 
    &\bb{4.2} & \bb{5.7}& \bb{6.6}
    &  \bb{4.2} & \bb{5.7} & \bb{6.6}& \bb{0.0549} \\
\hline[2pt]
\end{tblr}
}

\end{table}
 
\begin{table}[!htbp]
\centering
\caption{\small \textbf{Ill-Conditioned Tensor or Positive and Directional Tensors: }  \texttt{a} $=80$, \texttt{b} $=80$, and $\sigma =80$.  In the \textbf{first table}, $T$ is a diagonal ill-conditioned tensor with diagonal entries uniformly distributed in $[10^{-10}, 10^3]$. In the \textbf{second table}, we consider positive and directional tensors, where $T$ is a diagonal tensor with entries that are uniformly distributed in the range of $[0,40]$.   }
     \label{table large tensor}
{\scriptsize
\begin{tblr}{
  colspec = {c|rrr|rrr|rrr|r},
}
\hline[2pt]
   & \SetCell[c=3]{}Relative Min. & & & \SetCell[c=3]{}Iterations. & & &  \SetCell[c=3]{}Evaluations.& & & CPU time \\
\hline[1pt]
Method  &   $n=5$ &  50 &  100 & 5 &  50 &  100  &   5 &  50 &  100 & (s)\\
\hline[0.5pt]
    \textbf{\color{red}ARC}  & -1 & -1 & -1 
    &  9 & 11.7 & 15.1
    &  7.9 & 9.6 & 11.6 & 0.077 \\
   \textbf{\color{blue}QQR-v1}  &-1 & -1 & -1 
   & 13.3 & 11.8 & 11.6
   & 12.5 & 11.5 & 11.5 &0.090 \\
   \textbf{\color{blue}QQR-v2} & -1 & -1 & -1 
    &  \bb{7.6}  & \bb{8.9} & \bb{8.4}
    &  \bb{7.3} & \bb{8.5} & \bb{8.4} & \bb{0.071} \\
\hline[0.5pt]
\hline[0.5pt]
    \textbf{\color{red}ARC}  &  -1 & -1 & -1 
    & 7.1 & 8.5 & 14
    &  6.2 & 6.5 & 8 & 0.0501 \\
   \textbf{\color{blue}QQR-v1}   &  -1 & -1 & -1 
   & 6.9  & 8 & 8.1
   & 6.8 & 7.5 & 7.5 & 0.0509\\
   \textbf{\color{blue}QQR-v2} &  -1 & -1 & -1 
    &\bb{4.1} & \bb{4.4}& \bb{4}
    &  \bb{4.1} & \bb{4.4} & \bb{4}& \bb{0.0293} \\
\hline[2pt]
\end{tblr}
}
\end{table}

\newpage
\clearpage

\subsection{Minimizing General Objectives Using the AR$3$ algorithm}

This section provides a preliminary benchmark comparison of second and third-order methods on  20 problems from the unconstrained minimization problems in the Moré, Garbow, and Hillstrom (MGH) test set \cite{more1981testing}. The problem initialization, as well as the first, second, and third-order derivatives, were explicitly calculated and detailed in \cite{birgin2018fortran}\footnote{The original implementation for the test problems was written in Fortran, but we re-coded the test problems in Matlab using a Fortran-Matlab converter.}. All tests were conducted on an Intel(R) Core(TM) i7-4770 CPU @ 3.40 GHz processor with 16 GB of RAM. 

The performance of the following methods was compared:\textit{ the Adaptive Regularization with the Cubic method} (referred to as \texttt{ARC}), \textit{the AR\(3\) method with the inner subproblem solver ARC} (referred to as \texttt{AR3 with ARC}), and \textit{the AR\(3\) method with the inner subproblem solver QQR} (referred to as \texttt{AR3 with QQR}). 
The ARC algorithm is implemented as described in \cite{cartis2011adaptive}\footnote{URL: \text{https://www.galahad.rl.ac.uk/packages/}}. For third-order methods, numerical experiments were carried out to minimize a general objective function \( f(x) \) using the AR\(3\) algorithm as specified in \cite{cartis2020concise}, with the subproblem defined as in \eqref{subprob}. The subproblem is solved using either the ARC method or the \texttt{QQR-v2} method.
In Appendix \ref{Appendix numerical}, Tables \ref{full model 1}--\ref{full model 3}, we provide numerical details for each problem, where \texttt{eval.} and \texttt{iter.} denote the number of successful outer iterations and the total number of outer iterations, respectively; \texttt{inner eval.} represents the average number of successful inner iterations required to solve each subproblem, while \texttt{inner iter.} represents the average total number of inner iterations required to solve each subproblem. The performance ratio parameters for AR$3$, ARC, and QQR are set as follows: $\rho_1 = 0.1$, $\rho_2 = 0.9$, $\eta_0 = 0.5$, and $\eta_1 = 2$. The rest of the parameters in QQR are set as in Section \ref{sec: algo for QQRv2}; the tolerances and maximum iterations for the subproblem solver and outer iterations are set to $10^{-6}$ and $10^{-3}$, and 1000 and 3000, respectively.

Performance profile plots are a standard tool for comparing the performance of different algorithms on a collection of test problems using various metrics \cite{dolan2002benchmarking, dolan2002benchmarking, gould2016note}. In this work, we follow the approach and notation in \cite{birgin2017use}, letting \(\Gamma_{\hat{m}}(\tau_m)\) denote the proportion of problems for which the respective method (ARC, AR3 + ARC or AR3 + QQR)  was able to find an $\epsilon_g$-approximate solution within a given tolerance, with a `cost' up to \(\tau_m\) times the `cost' required by the most efficient method; here, the `cost' will be measured in terms of iteration, evaluation counts and CPU time. Specifically, if a method $\hat{m}$ finds an $\epsilon_g$--approximation to a solution of a problem $\hat{j}$, we define the cost of finding such an approximate solution as $t_{\hat{m} \hat{j}}(\epsilon_g)$. 
If the problem is not solved within the maximum iteration limit, we set $t_{\hat{m} \hat{j}}(\epsilon_g) = +\infty$. Note that here, $t_{\hat{m} \hat{j}}$  represents CPU time, number of iterations, and function- and derivative-evaluations. 
Define 
$
t^{\hat{j}}_{\min} (\epsilon_g) = \min_{\hat{m} \in \{1,2,3\}} t_{\hat{m} \hat{j}},
$
where $\hat{m} \in \{1,2,3\}$ represents the three methods: \text{ARC, AR3 + ARC, and AR3 + QQR}, respectively. Note that we have $t^{\hat{j}}_{\min} (\epsilon_g) = +\infty$ if none of the evaluated methods reach desired accuracy within allowed iteration limit.
The efficiency and robustness of a method can be evaluated by analyzing the curve
\[
\Gamma_{\hat{m}}(\tau_m) := \frac{1}{\tilde{q}} \# \bigg\{ \hat{j} \in \{1, \dotsc, \tilde{q}\} \ \big| \ t_{\hat{m} \hat{j}} (\epsilon_g) \le \tau_m t^{\hat{j}}_{\min} (\epsilon_g)\bigg\},
\]
where $\#$ denotes the cardinality of the set and $\tilde{q}$ represents the total number of problems.

\begin{figure}[ht!]
  \begin{center}
    \includegraphics[width=16cm]{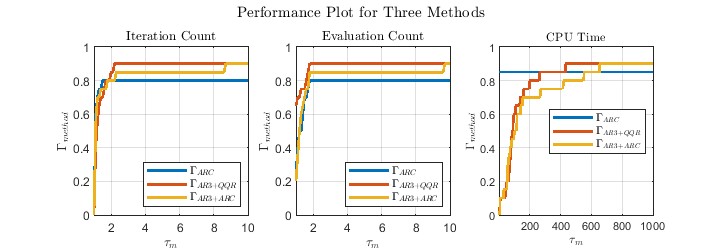}
\caption{\small Performance profile plots for the three methods. $\tau_m \in [1,10]$ is used for iteration and evaluation counts, while $\tau_m \in [1,10^3]$ is used for CPU time. A zoomed-in plot for $\tau_m \in [1,3]$ for iteration and evaluation counts, and $\tau_m \in [1,100]$ for CPU time, is provided in Appendix \ref{Appendix numerical}, Figure \ref{fig performance plot appendix}.}
    \label{fig performance plot long}
  \end{center}
\end{figure}

From Figure \ref{fig performance plot long}, we observe that, in general, third-order methods require at most as many outer function and derivative evaluations as second-order methods, and sometimes less. Between the two third-order methods, \texttt{AR3 with QQR} requires the least number of function and derivative evaluations. On average, across all problems, the \texttt{ARC} method requires 628 function and derivative evaluations, while \texttt{AR3 with ARC} requires 443 on average, and \texttt{AR3 with QQR} requires only 275 on average. 
The most significant differences in the number of iterations and evaluations typically occur for problems with a dominant tensor term. For instance,  for Problem 4 (`Brown badly scaled problem') in the Moré, Garbow, and Hillstrom’s test set \cite{more1981testing}, the tensor term has a dominant single entry of size  $10^8$, while all other entries are relatively small, typically of order $1$. For this problem, \texttt{AR3 with ARC} requires 1372 evaluations, whereas \texttt{AR3 with QQR} requires only 31 evaluations. 
Another example is Problem 19, where \texttt{AR3 with QQR} takes 146 evaluations compared to 800 evaluations required by \texttt{AR3 with ARC}. The primary reason for this difference is that the local minimum of the subproblem found by the QQR solver provides a greater descent in the objective function. Out of all the problems,  Problem 10 and  Problem 20 are not solved by any of the three methods within the maximum limit of 3000 iterations; this behaviour of Problem 10 is also reported in \cite{birgin2017use, Cartis2024Efficient}.

Lastly, we note that in terms of   CPU time, the performance of \texttt{AR3 + QQR} and \texttt{AR3 + ARC} is comparable. The primary reason is that \texttt{AR3 + QQR} requires fewer iterations and evaluations  per subproblem compared to \texttt{AR3 + ARC}, with \texttt{AR3 + QQR} averaging 1.92 iterations and 1.80 evaluations, while \texttt{AR3 + ARC} averages 6.47 iterations and 4.08 evaluations.
When compared to second-order methods, the \texttt{AR3} method generally requires more CPU time than the ARC method.
Currently, both the quadratic-quartic model and the ARC model are solved using Newton's method applied to a secular equation, which then uses Cholesky factorization of the ensuing linear system.
 Other alternatives to the latter are possible, such as computing, for small-scale problems, the complete eigenvalue-eigenvector decomposition. 
Similarly to existing approaches for trust-region and regularization subproblems (\cite[Ch. 6]{cartis2011adaptive} and \cite[Ch. 8, 10]{cartis2021scalable}), scalable iterative algorithms could be designed here, based on Krylov methods, eigenvalue formulations, or subspace optimization. A detailed analysis of these alternative approaches is deferred to future work.

The numerical results in this section only provide a preliminary numerical illustration of applying AR$3$ methods to objective functions under the simplest setup for adjusting the regularization parameter  \cite{cartis2020concise}, paving the way for future studies.  
Further work is needed to develop a competitive, optimized implementation of  AR$3$ methods, which we are undertaking in a dedicated way in  \cite{Cartis2024Efficient}.

\section{Conclusions and Future Work}
\label{Conclusion}

In conclusion, this paper presented the QQR methods, second-order methods that aim to achieve high-accuracy minimization of nonconvex quartically-regularised cubic polynomials, specifically the AR$3$ sub-problem. The QQR method is an iterative algorithm that operates by minimizing a sequence of quadratic polynomials with quartic regularization at the set of points $\{s^{(i)}\}_{i \ge 0}$.  This method can be viewed as a nonconvex extension of Nesterov's method proposed in \cite{Nesterov2022quartic}, which was originally designed for convex $m_3$ functions. Notably, when applied to convex $m_3$ functions, the QQR method exhibits the same linear convergence profile as Nesterov's method. For locally convex iterations, the QQR method achieves either a linear rate of convergence in gradient norm or an error reduction of at least $\mathcal{O}(\epsilon^{4/3})$ per iteration. In the case of locally nonconvex iterations, the QQR method provides an error reduction of at least $\mathcal{O}(\epsilon)$. 
The conducted numerical experiments validate the effectiveness of the QQR methods, as they perform competitively with state-of-the-art approaches such as ARC. In many instances, the QQR methods require fewer iterations or function evaluations and can find lower minima. Additionally, the QQR methods can handle both convex and nonconvex quartically-regularised cubic polynomials while effectively addressing challenges related to ill-conditioning in the Hessian term and the tensor term. Previous numerical studies \cite{birgin2017use, chow1989derivative, schnabel1971tensor, schnabel1991tensor} have indicated that the AR$3$ method provides enhanced optimization efficiency in terms of function and derivative evaluations when compared to first- or second-order methods. However, these studies have primarily focused on utilizing conventional unconstrained minimization methods for the optimization of AR$3$ sub-problems. On the contrary, the QQR methods are designed specifically for the AR$3$ sub-problem and approximate the third-order tensor term using a linear combination of quadratic and quartic terms. 

Developing efficient methods for minimizing polynomial models with tensor terms is crucial for high-order methods. The proposed QQR method provides a way to approximate and bound the third-order tensor term with a quadratic model with quartic regularization. By proposing efficient AR$3$ minimization algorithms, this paper brings high-order tensor methods closer to practical use.

\smallbreak 
\noindent
\textbf{Acknowledgments}: This work was supported by the Hong Kong Innovation and Technology Commission (InnoHK Project CIMDA).

\scriptsize{
\bibliographystyle{plain}
\bibliography{sample.bib}
}

\normalsize
\newpage

\appendix

\section{Numerical Results of the Experiment with MGH Test Set}
\label{Appendix numerical}
\small
\begin{table}[h!]
\centering
\caption{\textbf{Optimization Results using ARC for the MGH test set.} The parameters are $\rho_1 = 0.1$, $\rho_2 = 0.9$, $\eta_0 = 0.5$, and $\eta_1 = 2$ for both the inner and outer solvers. The tolerances for the inner and outer layers are set to $10^{-6}$ and $10^{-3}$, respectively. The symbol \^ indicates that the local minimum was not determined within the 3000 maximum iterations but can be attained if the iteration limit is increased to 5000. The symbol $*$ indicates that the local minimum was not determined even within 5000 maximum iterations.}
\begin{tabular}{@{}ccccc@{}}
\toprule
\textbf{No.} & \textbf{ARC Iter.} & \textbf{ARC Eval.} & \textbf{Func. Min.} & \textbf{Time (s)} \\ 
\midrule
1  & 2317 & 2317 & 0.00000      & 0.167  \\
2  & 70   & 70   & 48.98425     & 0.007  \\
3  & 26   & 13   & 0.96031      & 0.006  \\
4  & 1372 & 1372 & 0.00000      & 0.102  \\
5  & 57   & 57   & 0.45389      & 0.007  \\
6  & 104  & 104  & 124.3622     & 0.011  \\
7  & 779  & 778  & 0.00000      & 0.059  \\
8  & 707  & 707  & 0.00830      & 0.055  \\
9  & 24   & 17   & 0.00000      & 0.004  \\
10* & -   & -    & -            & -  \\
11 & 29   & 28   & 0.00419      & 0.007  \\
12 & 24   & 24   & 0.00000      & 0.003  \\
13 & 209  & 209  & 0.00000      & 0.017  \\
14 & 53   & 40   & 0.00000      & 0.008  \\
15 & 24   & 22   & 0.00034      & 0.003  \\
16 & 91   & 91   & 85,822       & 0.009  \\
17\^ & 3000 & 2988 & 0.00007     & 0.232  \\
18 & 13   & 13   & 0.30637      & 0.003  \\
19 & 96   & 83   & 0.00449      & 0.055  \\
20\^ & 3000 & 3001 & 0.00201     & 0.251  \\
\midrule
Total & \textbf{631} & \textbf{628} &             & \textbf{0.070 } \\
\bottomrule
\end{tabular}

\label{full model 1}
\end{table}

\begin{table}[h!]
\centering
\caption{\textbf{Optimization Results using AR$3$ with QQR as subproblem solver.} Parameters are the same as Table \ref{full model 1}.}
\begin{tabular}{ccccccc}
\toprule
\textbf{No.} & \textbf{AR$3$ Iter.} & \textbf{AR$3$ Eval.} & \textbf{Inner Iter.} & \textbf{Inner Eval.} & \textbf{Func. Min.} & \textbf{Time} \\ 
\midrule
1   & 2335 & 1754 & 1.39 & 1.32 & 0.00000 & 19.10 \\
2   & 70   & 49   & 1.47 & 1.47 & 48.98425 & 0.55  \\
3   & 31   & 12   & 1.90 & 1.90 & 0.96031 & 0.29  \\
4   & 32   & 31   & 1.94 & 1.94 & 0.00000 & 0.29  \\
5   & 71   & 51   & 1.21 & 1.21 & 0.45352 & 0.47  \\
6   & 118  & 87   & 3.22 & 3.05 & 124.3622 & 1.48  \\
7   & 694  & 568  & 1.12 & 1.12 & 0.00000 & 4.56  \\
8   & 794  & 417  & 1.01 & 1.01 & 0.00829 & 4.75  \\
9   & 27   & 13   & 1.19 & 1.19 & 0.00000 & 0.18  \\
10* & -    & -    & -    & -    & -       & -     \\
11  & 57   & 47   & 1.11 & 1.11 & 0.00418 & 0.36  \\
12  & 25   & 24   & 1.28 & 1.28 & 0.00000 & 0.36  \\
13  & 186  & 119  & 1.16 & 1.15 & 0.00001 & 1.3   \\
14  & 114  & 62   & 7.78 & 6.63 & 0.00000 & 2.91  \\
15  & 22   & 21   & 1.05 & 1.05 & 0.00036 & 0.64  \\
16  & 92   & 85   & 2.59 & 2.49 & 85,822  & 0.97  \\
17  & 174  & 147  & 3.44 & 2.70 & 0.00007 & 2.06  \\
18  & 17   & 16   & 1.53 & 1.53 & 0.30637 & 0.17  \\
19  & 159  & 146  & 1.04 & 1.04 & 0.00448 & 1.19  \\
20\^ & 3001 & 1581 & 1.00 & 1.00 & 0.00220 & 86.29 \\
\midrule
\textbf{Total} & \textbf{422}  & \textbf{275}  & \textbf{1.92} & \textbf{1.80} & - & \textbf{6.74} \\
\bottomrule
\end{tabular}

\label{full model 2}
\end{table}

\begin{table}[h!]
\centering
\caption{\textbf{Optimization Results using AR$3$ with ARC as subproblem solver}. Parameters are the same as Table \ref{full model 1}.}
\begin{tabular}{ccccccc}
\toprule
\textbf{No.} & \textbf{AR$3$ Iter.} & \textbf{AR$3$ Eval.} & \textbf{Inner Iter.} & \textbf{Inner Eval.} & \textbf{Func. Min.} & \textbf{Time} \\ 
\midrule
1   & 2354 & 2311 & 4.16 & 3.79 & 0.00000     & 22.03   \\
2   & 58   & 57   & 2.74 & 2.74 & 48.98425    & 0.39    \\
3   & 31   & 12   & 8.55 & 4.84 & 0.96031     & 0.45    \\
4   & 32   & 31   & 10.69 & 5.75 & 0.00000    & 0.57    \\
5   & 58   & 57   & 2.45 & 2.45 & 0.45390     & 0.33    \\
6   & 119  & 99   & 17.24 & 15.60 & 124.36218 & 4.23    \\
7   & 717  & 716  & 2.09 & 2.08 & 0.00000     & 3.45    \\
8   & 726  & 725  & 2.01 & 2.01 & 0.00829     & 3.41    \\
9   & 17   & 12   & 3.06 & 3.06 & 0.00000     & 0.13    \\
10* & -    & -    & -    & -    & -           & -       \\
11  & 50   & 40   & 13.50 & 4.14 & 0.00418    & 0.90    \\
12  & 26   & 25   & 2.50 & 2.50 & 0.00000     & 0.35    \\
13  & 208  & 206  & 2.17 & 2.15 & 0.00000     & 1.45    \\
14  & 119  & 63   & 13.13 & 7.59 & 0.00000    & 2.98    \\
15  & 23   & 22   & 2.83 & 2.30 & 0.00036     & 0.75    \\
16  & 92   & 85   & 3.53 & 3.15 & 85,822.20   & 0.79    \\
17  & 173  & 148  & 7.01 & 3.60 & 0.00007     & 2.28    \\
18  & 17   & 16   & 3.47 & 3.41 & 0.30637     & 0.20    \\
19  & 831  & 800  & 19.80 & 4.35 & 0.00070    & 21.82   \\
20\^  & 3001 & 3001 & 2.01 & 2.01 & 0.00201     & 142.33 \\
\midrule
\textbf{Total} & \textbf{455} & \textbf{443} & \textbf{6.47} & \textbf{4.08} & - & \textbf{10.99} \\
\bottomrule
\end{tabular}
\label{full model 3}
\end{table}

\begin{figure}[ht!]
  \begin{center}
    \includegraphics[width=16cm]{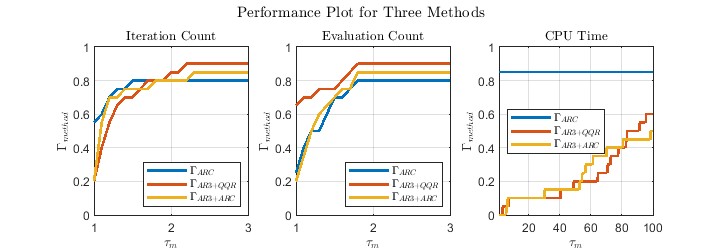}
    \caption{\small  Performance profile plots for three methods; $\tau_m \in [1,3]$ for iteration and evaluation counts and $\tau_m \in [1,100]$ for CPU time.}
    \label{fig performance plot appendix}
  \end{center}
\end{figure}

\end{document}